\documentclass[dvips,titlepage,letterpaper,12pt,fleqn]{article}
\usepackage[total={8.5in,11in}, top=1.0in, left=1.0in, right=1.0in, bottom=1.0in]{geometry}
\usepackage{graphicx}
\usepackage{amsmath}
\usepackage{amssymb}
\usepackage{amsfonts}
\usepackage{latexsym}
\usepackage{rotate}
\usepackage{float}
\usepackage{supertabular}
\usepackage{rotating}
\usepackage[frame]{xy}
\usepackage{fancybox}
\usepackage{tabularx}
\usepackage{multirow}
\usepackage{color}
\usepackage{mathrsfs}
\usepackage{epsfig}
\usepackage{lscape}
\usepackage{appendix}
\usepackage{tabu}
\usepackage{mathdots}
\usepackage{array,booktabs,calc}
\usepackage{adjustbox}
\usepackage[natbibapa]{apacite}
\usepackage{psfrag}
\usepackage{amsbsy}
\usepackage{amsthm}
\usepackage{epstopdf}
\usepackage{psfrag}
\usepackage{bigstrut}
\newcommand{\bbe}{\mbox{\boldmath$e$}}

\newcommand{\bA}{\mbox{\boldmath$A$}}

\newcommand{\bE}{\mbox{\boldmath$E$}}

\newcommand{\bH}{\mbox{\boldmath$H$}}

\newcommand{\bJ}{\mbox{\boldmath$J$}}
\newcommand{\bK}{\mbox{\boldmath$K$}}

\newcommand{\bP}{\mbox{\boldmath$P$}}
\newcommand{\bQ}{\mbox{\boldmath$Q$}}
\newcommand{\bR}{\mbox{\boldmath$R$}}

\newcommand{\bT}{\mbox{\boldmath$T$}}

\newcommand{\bV}{\mbox{\boldmath$V$}}
\newcommand{\bW}{\mbox{\boldmath$W$}}
\newcommand{\bX}{\mbox{\boldmath$X$}}
\newcommand{\bY}{\mbox{\boldmath$Y$}}
\newcommand{\bZ}{\mbox{\boldmath$Z$}}

\newcommand{\be}{\begin{equation*}}
\newcommand{\ee}{\end{equation*}}
\newcommand{\ben}{\begin{equation}}
\newcommand{\een}{\end{equation}}
\newcommand{\bea}{\begin{eqnarray*}}
\newcommand{\eea}{\end{eqnarray*}}
\newcommand{\bean}{\begin{eqnarray}}
\newcommand{\eean}{\end{eqnarray}}
\newcommand{\ba}{\begin{array}}
\newcommand{\ea}{\end{array}}
\newcommand{\lba}{\left[ \begin{array}}
\newcommand{\ear}{\end{array} \right]}
\newcommand{\vect}{\mbox{\rm vec}}

\newcommand{\cE}{{\rm I\hspace{-0.7mm}E}}
\newcommand{\cR}{{\rm I\hspace{-0.7mm}R}}

\newtheorem{definition}{Definition}

\newlength{\LL} 
\newlength{\NL} 
\newlength{\MD} 
\title{Computation of the State Bias and Initial States \\ for Stochastic State Space Systems in the General 2-D Roesser Model Form}
\author{Jos\'{e} A. Ramos \\ Nova Southeastern University \\
College of Engineering and Computing \\
Department of Engineering and Technology \\
3301 College Avenue \\
Fort Lauderdale, FL 33314 \\
Email: jr1284@nova.edu \\ and \\ Guillaume Merc\`{e}re \\
Universit\'{e} de Poitiers \\
Laboratoire d'Informatique et d'Automatique pour les Syst{\`e}mes \\
2 rue Pierre Brousse, b{\^a}timent B25, TSA 41105 \\
86073 Poitiers cedex 9, France \\
Email: guillaume.mercere@univ-poitiers.fr}

\begin{document}

\maketitle

\begin{abstract}
Recently \cite{Ramos2017a} presented a subspace system identification algorithm for
2-D purely stochastic state space models in the general Roesser form. However, since the exact problem requires an
oblique projection of $\bY_f^h$ projected onto $\bW_p^h$ along $\widehat{\bX}_f^{vh}$, where $\bW_p^h=
\lba{c}\widehat{\bX}_p^{vh} \\ \bY_p^h \ear$, this presents a problem since $\{\widehat{\bX}_p^{vh},\widehat{\bX}_f^{vh}\}$
are unknown. In the above mentioned paper, the authors found that by doing an orthogonal projection $\bY_f^h/\bY_p^h$,
one can identify the future horizontal state matrix $\widehat{\bX}_f^{h}$ with a small bias due to the initial conditions
that depend on $\{\widehat{\bX}_p^{vh},\widehat{\bX}_f^{vh}\}$. Nevertheless, the results on modeling 2-D images were
very good despite lack of knowledge of $\{\widehat{\bX}_p^{vh},\widehat{\bX}_f^{vh}\}$. In this note we delve into the
bias term and prove that it is insignificant, provided $i$ is chosen large enough and the vertical and horizontal states are uncorrelated. 
That is, the cross covariance of the state estimates $x_{r,s}^{h}$ and $x_{r,s}^{v}$ is zero, or $P_{hv}=0_{n_x\times n_x}$ and $P_{vh}=0_{n_x\times n_x}$. Our simulations use $i=30$. We also present
a second iteration to improve the state estimates by including the vertical states computed from a vertical data processing step,
i.e., by doing an orthogonal projection $\bY_f^v/\bY_p^v$. In this revised algorithm we include a step to compute the initial
states. This new portion, in addition to the algorithm presented in \cite{Ramos2017a}, forms a complete 2-D stochastic subspace system identification algorithm.
\end{abstract}

\section{Problem Formulation}
The general 2-D stochastic Roesser model has the state-space form
\begin{subequations}
\bean \label{xh}
x^h_{r+1,s} & = & A_1x^h_{r,s} + A_2x^v_{r,s} + w^h_{r,s} \\ \label{xv}
x^v_{r,s+1} & = & A_3x^h_{r,s} + A_4 x^v_{r,s} + w^v_{r,s} \\ \label{y}
y_{r,s} & = & C_1x^h_{r,s} + C_2x^v_{r,s} + v_{r,s},
\eean
\end{subequations}
where $x^h_{r,s}\in\mathbb{R}^{n_{h}}$, $x^v_{r,s}\in\mathbb{R}^{n_{v}}$,
and $y_{r,s}\in\mathbb{R}^{n_y}$ denote, respectively, the local horizontal state, local vertical state,
and output vectors at the $(r,s)^{th}$ location of a finite domain $\mathbb{D}=\{(r,s)\;|\; 0\le r\le N\; \mbox{and}\; 0\le s\le M\}$.
The system matrices $\{A,C\}$, given by
\begin{subequations}
\bean
A & = & \lba{c|c} A_1 & A_2 \\ \hline A_3 & A_4 \ear,\;\;\;
C \;=\; \lba{c|c} C_1 & C_2 \ear,
\eean
\end{subequations}
have partitioned dimensions $A_1\in\mathbb{R}^{n_h\times n_h}$, $A_2\in\mathbb{R}^{n_h\times n_v}$,
$A_3\in\mathbb{R}^{n_v\times n_h}$, $A_4\in\mathbb{R}^{n_v\times n_v}$, $C_1\in\mathbb{R}^{n_y \times n_h}$, and $C_2\in\mathbb{R}^{n_y\times n_v}$.
The noise vectors $w^h_{r,s}\in\mathbb{R}^{n_h}$, $w^v_{r,s}\in\mathbb{R}^{n_v}$, and $v_{r,s}\in\mathbb{R}^{n_y}$
are assumed to be white Gaussian noise processes with mean and joint covariance matrix given, respectively, by
\renewcommand*{\arraystretch}{1.25}
\begin{subequations}
\bean
\cE\left\{\lba{c} w^h_{r,s} \\ \hline w^v_{r,s} \\ \hline v_{r,s} \ear \right\} & = & \lba{c} 0_{n_h \times 1}
\\ \hline 0_{n_v \times 1} \\ \hline 0_{n_y \times 1} \ear \\
\cE\left\{\lba{c} w^h_{r,s} \\ w^v_{r,s} \\ \hline v_{r,s} \ear \lba{cc|c} \left(w^h_{r',s'}\right)^{\top} &
\left(w^v_{r',s'}\right)^{\top} & v^{\top}_{r',s'} \ear \right\}
& = & \lba{cc|c} Q_{hh} & Q_{hv} & S_h \\ Q_{vh} & Q_{vv} & S_v \\ \hline S^{\top}_h & S^{\top}_v & R \ear
\cdot \delta_{r-r'}\cdot \delta_{s-s'} \nonumber \\
& = & \lba{c|c} Q & S \\ \hline S^{\top} & R \ear \cdot \delta_{r-r'}\cdot \delta_{s-s'},
\eean
\end{subequations}
\noindent where $Q_{hh}\in\mathbb{R}^{n_h\times n_h}$, $Q_{hv}\in\mathbb{R}^{n_h\times n_v}$, $Q_{vh}\in\mathbb{R}^{n_v\times n_h}$,
$Q_{vv}\in\mathbb{R}^{n_v\times n_v}$, $S_{h}\in\mathbb{R}^{n_h\times n_y}$, $S_{v}\in\mathbb{R}^{n_v\times n_y}$,
and $R\in\mathbb{R}^{n_y\times n_y}$, $n_x=n_h+n_v$ is the dimension of the combined system, $\cE$ is the expectation operator, $M^{\top}$ denotes the transpose of $M$, $\delta_{k-k'}$ is the Kronecker delta function, $0_{m \times n}$ denotes an $(m \times n)$ matrix with all its elements equal to zero, and $\{Q,R,S\}$ are the covariance and cross-covariance matrices of the noise terms.

The noise and state vectors are uncorrelated with each other, i.e.,
\renewcommand*{\arraystretch}{1.00}
\begin{subequations}
\bean
\cE\left\{x^h_{r,s}\lba{c|c|c} \left(w^h_{r',s'}\right)^{\top} & \left(w^v_{r',s'}\right)^{\top} & v_{r',s'}^{\top} \ear \right\}
& = & 0_{n_h \times (n_x+n_y)}, \; \forall \; r' \ge r \; \mbox{and} \; s' \ge s \\
\cE\left\{x^v_{r,s}\lba{c|c|c} \left(w^h_{r',s'}\right)^{\top} & \left(w^v_{r',s'}\right)^{\top} & v_{r',s'}^{\top} \ear \right\}
& = & 0_{n_v \times (n_x+n_y)}, \; \forall \; r' \ge r \; \mbox{and} \; s' \ge s.
\eean
\end{subequations}
Furthermore, the states $x^h_{r,s}$ and $x^v_{r,s}$ evolve with the following statistical properties: zero mean
\renewcommand*{\arraystretch}{1.25}
\bean
\cE\left\{\lba{c} x^h_{r,s} \\ \hline x^v_{r,s} \ear \right\} & = & \lba{c} 0_{n_h \times 1} \\ \hline 0_{n_v \times 1} \ear,
\;\; r=0,1,\ldots,N \; \mbox{and} \; s=0,1,\ldots,M
\eean
and positive definite state covariance matrix
\bean \label{Pih}
\Pi & = & \cE\left\{\lba{c} x^h_{r,s} \\ \hline x^v_{r,s} \ear \lba{c|c} (x^h_{r,s})^{\top} & (x^v_{r,s})^{\top} \ear \right\} \;=\;
\lba{c|c} \Pi_h & 0_{n_h \times n_v} \\ \hline 0_{n_v \times n_h} & \Pi_v \ear.
\eean
Let us now define the covariance of the state update as
\bean
\Pi' & = & \cE\left\{\lba{c} x^h_{r+1,s} \\ \hline x^v_{r,s+1} \ear \lba{c|c} \left(x^h_{r+1,s}\right)^{\top} & \left(x^v_{r,s+1}\right)^{\top} \ear \right\} \;=\;
\lba{c|c} \Pi_{h} & \Pi_{hv} \\ \hline \Pi^{\top}_{hv} & \Pi_{v} \ear,
\eean
where $\Pi_{hv}=A_1\Pi_{h}A_3^{\top} + A_2\Pi_{v}A_4^{\top} + Q_{hv}$, $\Pi_{vh}=\Pi^{\top}_{hv}$, and the dimensions are
$\Pi_{h}\in\mathbb{R}^{n_h\times n_h}$, $\Pi_{hv}\in\mathbb{R}^{n_h\times n_v}$, $\Pi_{vh}\in\mathbb{R}^{n_v\times n_h}$, and
$\Pi_{v}\in\mathbb{R}^{n_v\times n_v}$. The state covariance update equation becomes
\renewcommand*{\arraystretch}{1.00}
\bean \label{PP1}
\Pi' & = & A\Pi A^{\top} + Q,
\eean
where $\Pi=\Pi^{\top}$ and $\Pi'=(\Pi')^{\top}$.
Note that \eqref{PP1} is not a matrix Lyapunov state covariance equation since $\Pi'\ne \Pi$. However, by partitioning
\eqref{PP1}, one can decompose it into a pair of coupled horizontal and vertical matrix Lyapunov type equations \citep{Ramos2016b}. Nevertheless, one can enforce the constraint $\Pi_{hv}=0_{n_h \times n_v}$, which results in the joint matrix Lyapunov equation
\bea
\lba{c|c} \Pi_{h} & 0_{n_h \times n_v} \\ \hline 0_{n_v \times n_h} & \Pi_{v} \ear & = & \lba{c|c} A_1 & A_2 \\ \hline A_3 & A_4 \ear
\lba{c|c} \Pi_{h} & 0_{n_h \times n_v} \\ \hline 0_{n_v \times n_h} & \Pi_{v} \ear \lba{c|c} A_1 & A_2 \\ \hline A_3 & A_4 \ear^{\top} + \lba{c|c} Q_{hh} & Q_{hv} \\ \hline
Q_{vh} & Q_{vv} \ear,
\eea
or, more compactly,
\bean \label{Pi}
\Pi & = & A\Pi A^{\top} + Q.
\eean
Throughout the rest of this note we will use the symbol $>0$ $(\ge 0)$ to indicate that a matrix is positive definite $($positive semi-definite$)$. Model \eqref{xh} -- \eqref{y} then satisfies the following constraints, also known as the positive real conditions:
\renewcommand*{\arraystretch}{1.25}
\bean \label{PR1}
\lba{c|c} Q & S \\ \hline S^{\top} & R \ear & \ge & 0, \;\;\; Q \ge 0,\;\;\; R > 0,\;\;\;\Pi > 0.
\eean

The 2-D output autocovariance sequence $\Lambda_{k,m}\in\mathbb{R}^{n_y\times n_y}$ is given in terms of the Markov parameters
of the system as
\bean \label{Lambda}
\Lambda_{k,m} & = & \cE\left\{y_{r+k,s+m}\,y_{r,s}^{\top}\right\} \;=\;
\left\{\ba{ll} C_1\Pi_hC_1^{\top} + C_2\Pi_vC_2^{\top} + R, & \mbox{if}\;k=0,\;m=0 \\
C_1A_1^{k-1}G_1, & \mbox{if}\;k\ge 1,\;m=0 \\
C_2A_4^{m-1}G_2, & \mbox{if}\;k=0,\;m\ge 1 \\
CA^{k-1,m}G^{1,0} +
CA^{k,m-1}G^{0,1}, & \mbox{if}\;k\ge 1,\;m\ge 1, \ea \right.
\eean
where $G_1$ and $G_2$ are defined, respectively, as the horizontal and vertical partitions
of the matrix $G\in\mathbb{R}^{n_x \times n_y}$, obtained from
\begin{subequations}
\bean
G & = & \cE\left\{\lba{c} x^h_{r+1,s} \\ \hline x^v_{r,s+1} \ear y^{\top}_{r,s} \right\} \;=\; A\Pi C^{\top}+S,
\eean
with $G_1\in\mathbb{R}^{n_h \times n_y}$ and $G_2\in\mathbb{R}^{n_v \times n_y}$ given as
\bean \label{G1}
G_1 & = & A_1\Pi_h C_1^{\top}+A_2\Pi_v C_2^{\top}+S_h \\ \label{G2}
G_2 & = & A_3\Pi_h C_1^{\top}+A_4\Pi_v C_2^{\top}+S_v,
\eean
\end{subequations}
and
\bea
G^{1,0} & = & \lba{c} G_1 \\ \hline 0_{n_v \times n_y} \ear,\;\;\;
G^{0,1} \;=\; \lba{c} 0_{n_h \times n_y} \\ \hline G_2 \ear \\
A^{0,0} & = & I_{n},\;\;\;
A^{1,0} \;=\; \lba{c|c} A_1 & A_2 \\ \hline 0_{n_v \times n_h} & 0_{n_v \times n_v} \ear,\;\;\;\;\;
A^{0,1} \;=\; \lba{c|c} 0_{n_h \times n_h} & 0_{n_h \times n_v} \\ \hline A_3 & A_4 \ear \\
A & = & A^{1,0}+A^{0,1},\;\;\;\;\;
A^{k,m} \;=\; A^{1,0}A^{k-1,m}+A^{0,1}A^{k,m-1},\;\;\;\;\mbox{for}\;\;(k,m) > (0,0) \\
A^{-k,m} & = & A^{k,-m} \;=\; 0_{n_x\times n_x},\;\;\;\;\mbox{for}\;\;\;\;k \ge 1,\; m \ge 1.
\eea
The problem can now be stated as follows:
\begin{definition}
Given a data matrix $\bY\in\mathbb{R}^{n_y(N+1) \times (M+1)}$ corresponding to the output sequence $y_{r,s}\in\mathbb{R}^{n_y}$, for $r=0,1,\ldots,N$ and
$s=0,1,\ldots,M$, find: $(i)$ the system orders $n_h$ and $n_v$ such that $n_x=n_h+n_v$,
$(ii)$ parameter matrices
$\{A,C,G\}$ up to a similarity transformation, $(iii)$ covariance matrices $\{\Pi,Q,R,S\}$,
and $(iv)$ the initial conditions $\{x^h_{0,s}\}_{s=0}^M$ and
$\{x^v_{r,0}\}_{r=0}^N$, subject to the constraints \eqref{PR1}, so that the $2nd$-order
statistics of the output of the system match those of the given output data.
\end{definition}

In order to simplify the analysis, we also formulate the problem in the innovations form as
\renewcommand*{\arraystretch}{1.00}
\begin{subequations}
\bean \label{xh1}
{\widehat x}^{h}_{r+1,s} & = & A_1{\widehat x}^{h}_{r,s} + A_2{\widehat x}^{v}_{r,s} + K_{1}e_{r,s} \\ \label{xh2}
{\widehat x}^{v}_{r,s+1} & = & A_3{\widehat x}^{h}_{r,s} + A_4{\widehat x}^{v}_{r,s} + K_{2}e_{r,s} \\ \label{yh1}
y_{r,s} & = & C_1{\widehat x}^{h}_{r,s} + C_2{\widehat x}^{v}_{r,s} + e_{r,s},
\eean
\end{subequations}
where ${\widehat x}^{h}_{r,s}\in\mathbb{R}^{n_h}$ and ${\widehat x}^{v}_{r,s}\in\mathbb{R}^{n_v}$ are,
respectively, the horizontal and vertical state estimates, with
state estimate covariance matrices $P_{h}=\cE\left\{{\widehat x}^{h}_{r,s}\left({\widehat x}^{h}_{r,s}\right)^{\top}\right\}\in\mathbb{R}^{n_h\times n_h}$ and
$P_{v}=\cE\left\{{\widehat x}^{v}_{r,s}\left({\widehat x}^{v}_{r,s}\right)^{\top}\right\}\in\mathbb{R}^{n_v \times n_v}$. Furthermore, we assume that
$P_{hv}=\cE\left\{{\widehat x}^{h}_{r,s}\left({\widehat x}^{v}_{r,s}\right)^{\top}\right\} = 0_{n_h \times n_v}$. These state estimate covariance matrices
satisfy the joint Riccati equation
\bean \label{Ric1}
P & = & APA^{\top} + (G-APC^{\top})(\Lambda_{0,0} - CPC^{\top})^{-1}(G-APC^{\top})^{\top},
\eean
where
\renewcommand*{\arraystretch}{1.25}
\bean \label{Pmat}
P & = & \lba{c|c} P_{h} & 0_{n_h \times n_v} \\ \hline 0_{n_v \times n_h} & P_{v} \ear
\eean
is a positive definite matrix.
We further define the innovations covariance matrix $R_e=\cE\left\{e_{r,s}e_{r,s}^{\top}\right\}$ as
\renewcommand*{\arraystretch}{1.00}
\bean
R_e & = & \Lambda_{0,0} - C_1P_hC_1^{\top} - C_2P_vC_2^{\top}
\eean
and state estimate errors and state estimate error covariance matrices, respectively, as
\bea
{\tilde x}^{h}_{r,s} & = & x^{h}_{r,s} - {\widehat x}^{h}_{r,s} \in\mathbb{R}^{n_h} \\
{\tilde x}^{v}_{r,s} & = & x^{v}_{r,s} - {\widehat x}^{v}_{r,s} \in\mathbb{R}^{n_v} \\
\Sigma_{h} & = & \cE\left\{{\tilde x}^{h}_{r,s}\left({\tilde x}^{h}_{r,s}\right)^{\top}\right\} \;=\; \Pi_{h}-P_{h} \in\mathbb{R}^{n_h \times n_h} \\
\Sigma_{v} & = & \cE\left\{{\tilde x}^{v}_{r,s}\left({\tilde x}^{v}_{r,s}\right)^{\top}\right\} \;=\; \Pi_{v}-P_{v} \in\mathbb{R}^{n_v \times n_v}.
\eea
Then $\Sigma_{h}$ and $\Sigma_{v}$ satisfy the joint Riccati equation
\bean \label{Ric2}
\Sigma & = & A\Sigma A^{\top} + Q + (A\Sigma C^{\top}+S)(C\Sigma C^{\top}+R)^{-1}(A\Sigma C^{\top}+S)^{\top},
\eean
where
\renewcommand*{\arraystretch}{1.25}
\bea
\Sigma & = & \lba{c|c} \Sigma_{h} & 0_{n_h \times n_v} \\ \hline 0_{n_v \times n_h} & \Sigma_{v} \ear
\;=\; \lba{c|c} \Pi_{h} & 0_{n_h \times n_v} \\ \hline 0_{n_v \times n_h} & \Pi_{v} \ear -
\lba{c|c} P_{h} & 0_{n_h \times n_v} \\ \hline 0_{n_v \times n_h} & P_{v} \ear.
\eea
Finally, the Kalman gain matrix is given by either of the following two expressions
\renewcommand*{\arraystretch}{1.00}
\begin{subequations}
\bean \label{K1}
K & = & (G-APC^{\top})(\Lambda_{0,0} - CPC^{\top})^{-1} \in\mathbb{R}^{n_x \times n_y} \\ \label{K2}
K & = & (A\Sigma C^{\top}+S)(C\Sigma C^{\top}+R)^{-1} \in\mathbb{R}^{n_x \times n_y},
\eean
where
\renewcommand*{\arraystretch}{1.25}
\bea
K & = & \lba{c} K_{1} \\ \hline K_{2} \ear,
\eea
\end{subequations}
with dimensions $K_1\in\mathbb{R}^{n_h \times n_y}$ and $K_2\in\mathbb{R}^{n_v \times n_y}$.
\renewcommand*{\arraystretch}{1.00}
\section{Horizontal Data Processing}
\noindent Let the horizontal and vertical past and future state matrices for $k=0,1,\ldots,M$ and
$N=2i+j-2$ be defined as
\bean \label{Xph}
\widehat{X}_p^h(k) & \triangleq & \lba{c|c|c|c|c} \widehat{x}_{0,k}^{h} & \widehat{x}_{1,k}^{h} & \widehat{x}_{2,k}^{h} & \cdots & \widehat{x}_{j-1,k}^{h} \ear\in\cR^{n_h \times j} \\ \label{Xfh}
\widehat{X}_f^h(k) & \triangleq & \lba{c|c|c|c|c} \widehat{x}_{i,k}^{h} & \widehat{x}_{i+1,k}^{h} & \widehat{x}_{i+2,k}^{h} & \cdots & \widehat{x}_{i+j-1,k}^{h} \ear\in\cR^{n_h \times j} \\ \label{Xpv}
\widehat{X}_p^{vh}(k) & \triangleq & \lba{ccccc} \widehat{x}_{0,k}^{v} & \widehat{x}_{1,k}^{v} & \widehat{x}_{2,k}^{v} & \cdots & \widehat{x}_{j-1,k}^{v} \\
\widehat{x}_{1,k}^{v} & \widehat{x}_{2,k}^{v} & \widehat{x}_{3,k}^{v} & \cdots & \widehat{x}_{j,k}^{v} \\
\widehat{x}_{2,k}^{v} & \widehat{x}_{3,k}^{v} & \widehat{x}_{4,k}^{v} & \cdots & \widehat{x}_{j+1,k}^{v} \\ \vdots & \vdots & \vdots & \iddots & \vdots \\
\widehat{x}_{i-1,k}^{v} & \widehat{x}_{i,k}^{v} & \widehat{x}_{i+1,k}^{v} & \cdots & \widehat{x}_{i+j-2,k}^{v} \ear\in\cR^{n_vi \times j} \\ \label{Xfv}
\widehat{X}_f^{vh}(k) & \triangleq & \lba{ccccc} \widehat{x}_{i,k}^{v} & \widehat{x}_{i+1,k}^{v} & \widehat{x}_{i+2,k}^{v} & \cdots & \widehat{x}_{i+j-1,k}^{v} \\
\widehat{x}_{i+1,k}^{v} & \widehat{x}_{i+2,k}^{v} & \widehat{x}_{i+3,k}^{v} & \cdots & \widehat{x}_{i+j,k}^{v} \\
\widehat{x}_{i+2,k}^{v} & \widehat{x}_{i+3,k}^{v} & \widehat{x}_{i+4,k}^{v} & \cdots & \widehat{x}_{i+j+1,k}^{v} \\
\vdots & \vdots & \vdots & \iddots & \vdots \\
\widehat{x}_{2i-1,k}^{v} & \widehat{x}_{2i,k}^{v} & \widehat{x}_{2i+1,k}^{v} & \cdots & \widehat{x}_{2i+j-2,k}^{v} \ear\in\cR^{n_vi \times j},
\eean
where throughout the sequel, subscripts $p$ and $f$ denote \emph{past} and \emph{future}, respectively,
superscripts $h$ and $v$ denote \emph{horizontal} and \emph{vertical}, respectively, $vh$ denotes
\emph{vertical from horizontal data processing}, and $i$ and $j$ are
fixed integer constants such that $j \gg i$ and $n_y i \gg \mbox{max}\{n_h,n_v\}$. \\

\noindent Likewise, we define the horizontal past and future innovations and output data matrices for $k=0,1,\ldots,M$ and
$N=2i+j-2$ as follows:
\bean \label{Ep}
E^h_{p}(k) & \triangleq & \lba{ccccc} e_{0,k} & e_{1,k} & e_{2,k} & \cdots & e_{j-1,k} \\
e_{1,k} & e_{2,k} & e_{3,k} & \cdots & e_{j,k} \\
e_{2,k} & e_{3,k} & e_{4,k} & \cdots & e_{j+1,k} \\
\vdots & \vdots & \vdots & \iddots & \vdots \\
e_{i-1,k} & e_{i,k} & e_{i+1,k} & \cdots & e_{i+j-2,k} \ear\in\cR^{n_y i \times j} \\ \label{Ef}
E^h_{f}(k) & \triangleq & \lba{ccccc} e_{i,k} & e_{i+1,k} & e_{i+2,k} & \cdots & e_{i+j-1,k} \\
e_{i+1,k} & e_{i+2,k} & e_{i+3,k} & \cdots & e_{i+j,k} \\
e_{i+2,k} & e_{i+3,k} & e_{i+4,k} & \cdots & e_{i+j+1,k} \\
\vdots & \vdots & \vdots & \iddots & \vdots \\
e_{2i-1,k} & e_{2i,k} & e_{2i+1,k} & \cdots & e_{2i+j-2,k} \ear\in\cR^{n_y i \times j} \\ \label{Yp}
Y^h_{p}(k) & \triangleq & \lba{ccccc} y_{0,k} & y_{1,k} & y_{2,k} & \cdots & y_{j-1,k} \\
y_{1,k} & y_{2,k} & y_{3,k} & \cdots & y_{j,k} \\
y_{2,k} & y_{3,k} & y_{4,k} & \cdots & y_{j+1,k} \\
\vdots & \vdots & \vdots & \iddots & \vdots \\
y_{i-1,k} & y_{i,k} & y_{i+1,k} & \cdots & y_{i+j-2,k} \ear\in\cR^{n_y i \times j} \\ \label{Yf}
Y^h_{f}(k) & \triangleq & \lba{ccccc} y_{i,k} & y_{i+1,k} & y_{i+2,k} & \cdots & y_{i+j-1,k} \\
y_{i+1,k} & y_{i+2,k} & y_{i+3,k} & \cdots & y_{i+j,k} \\
\vdots & \vdots & \vdots & \iddots & \vdots \\
y_{2i-1,k} & y_{2i,k} & y_{2i+1,k} & \cdots & y_{2i+j-2,k} \ear\in\cR^{n_y i \times j}.
\eean
One can easily show that the following equations are satisfied for $k=0,1,\ldots,M$
\bean \label{Yp1}
Y_{p}^h(k) & = & \Gamma_i^h \widehat{X}_{p}^h(k) + \Gamma_i^{vh} \widehat{X}_{p}^{vh}(k) + K_i^h E_{p}^h(k) \\ \label{Yf1}
Y_{f}^h(k) & = & \Gamma_i^h \widehat{X}_{f}^h(k) + \Gamma_i^{vh} \widehat{X}_{f}^{vh}(k) + K_i^h E_{f}^h(k) \\ \label{Xf1}
\widehat{X}_{f}^h(k) & = & A_1^i \widehat{X}_{p}^h(k) + \Phi_i^{vh} \widehat{X}_{p}^{vh}(k) + {\cal L}^h_i E_{p}^h(k),
\eean
where $\{\Gamma_i^h,\Phi_i^{vh},{\cal L}^h_i\}$ and other related matrices are defined as follows:
\bean
\Gamma_i^h & \triangleq & \lba{c} C_1 \\ C_1A_1 \\ \vdots \\ C_1A_1^{i-1} \ear \in\cR^{n_yi \times n_h} \\
{\cal L}^h_i & \triangleq & \lba{c|c|c|c} A_1^{i-1}K_1 & A_1^{i-2}K_1 & \cdots &
K_1 \ear \in\cR^{n_h \times n_y i} \\
\Phi_i^{h} & \triangleq & \lba{c|c|c|c} A_1^{i-1} & A_1^{i-2} & \cdots & I_{n_h} \ear \in\cR^{n_h \times n_hi}
\eean
\bean
\Phi_i^{vh} & = & \lba{c|c|c|c} A_1^{i-1}A_2 & A_1^{i-2}A_2 & \cdots & A_2 \ear \in\cR^{n_h \times n_vi} \;=\;
\Phi_i^{h} \cdot (I_i \otimes A_2) \\
\Theta_i^{h} & \triangleq & \lba{c} I_{n_h} \\ A_1 \\ \vdots \\ A_1^{i-1} \ear \in\cR^{n_hi \times n_h},
\eean
and $I_{k}$ denotes a $(k \times k)$
identity matrix. Finally, we define the lower triangular block Toeplitz matrices $\{G_{A_1}^{h},\Gamma_i^{vh},K^h_i\}$ as
\bean
G_{A_1}^{h} & \triangleq & \lba{cccc} 0_{n_y \times n_v} & & & \\ I_{n_h} & 0_{n_y \times n_v} & & \\
\vdots & \vdots  & \ddots \\ A_1^{i-2} & A_1^{i-3} & \cdots &
0_{n_y \times n_v} \ear \in\cR^{n_hi \times n_hi} \\ \label{GT}
\Gamma_i^{vh} & \triangleq & \lba{cccc} C_2 & & & \\ C_1A_2 & C_2 & & \\
\vdots & \vdots  & \ddots & \\ C_1A_1^{i-2}A_2 & C_1A_1^{i-3}A_2 & \cdots &
C_2 \ear \in\cR^{n_y i \times n_vi} \nonumber \\
& = & (I_i \otimes C_1)G_{A_1}^{h}(I_i \otimes A_2) + (I_i \otimes C_2) \\
K^h_i & \triangleq & \lba{cccc} I_{n_y} & & & \\ C_1K_1 & I_{n_y} & & \\
\vdots & \vdots & \ddots & \\ C_1A_1^{i-2}K_1 & C_1A_1^{i-3}K_1 & \cdots &
I_{n_y} \ear \in\cR^{n_y i \times n_y i} \nonumber \\
& = & (I_i \otimes C_1)G_{A_1}^{h}(I_i \otimes K_1) + (I_i \otimes I_{n_y}).
\eean
For the purpose of horizontal data processing we will work with
the equivalent horizontal subsystem
\bean
\widehat{x}_{r+1,s}^{h} & = & A_1\widehat{x}_{r,s}^{h}+A_2\widehat{x}_{r,s}^{v}+K_1e_{r,s} \\
y_{r,s} & = & C_1\widehat{x}_{r,s}^{h}+C_2\widehat{x}_{r,s}^{v}+e_{r,s},
\eean
for $r=0,1,\ldots,N$ and $s=0,1,\ldots,M$. However, at this point we need to make the following notational simplification
${\bar n}_h \triangleq n_h(M+1)$,
$\bar{n}_y \triangleq n_y(M+1)$, and ${\bar \jmath} \triangleq j(M+1)$.
Then, by defining
\bean
\bY_p^h & \triangleq & \lba{c|c|c|c} Y_p^h(0) & Y_p^h(1) & \cdots & Y_p^h(M) \ear\in\cR^{n_y \times {\bar \jmath}} \\
\bY_f^h & \triangleq & \lba{c|c|c|c} Y_f^h(0) & Y_f^h(1) & \cdots & Y_f^h(M) \ear\in\cR^{n_y \times {\bar \jmath}} \\
\widehat{\bX}_p^h & \triangleq & \lba{c|c|c|c} \widehat{X}_p^h(0) & \widehat{X}_p^h(1) & \cdots & \widehat{X}_p^h(M) \ear\in\cR^{n_h \times {\bar \jmath}} \\
\widehat{\bX}_f^h & \triangleq & \lba{c|c|c|c} \widehat{X}_f^h(0) & \widehat{X}_f^h(1) & \cdots & \widehat{X}_f^h(M) \ear\in\cR^{n_h \times {\bar \jmath}} \\
\widehat{\bX}_p^{vh} & \triangleq & \lba{c|c|c|c} \widehat{X}_p^{vh}(0) & \widehat{X}_p^{vh}(1) & \cdots & \widehat{X}_p^{vh}(M) \ear\in\cR^{n_vi \times {\bar \jmath}}
\eean
\bean
\widehat{\bX}_f^{vh} & \triangleq & \lba{c|c|c|c} \widehat{X}_f^{vh}(0) & \widehat{X}_f^{vh}(1) & \cdots & \widehat{X}_f^{vh}(M) \ear\in\cR^{n_vi \times {\bar \jmath}} \\
\bE_p^h & \triangleq & \lba{c|c|c|c} E_p^h(0) & E_p^h(1) & \cdots & E_p^h(M) \ear\in\cR^{n_y \times {\bar \jmath}} \\
\bE_f^h & \triangleq & \lba{c|c|c|c} E_f^h(0) & E_f^h(1) & \cdots & E_f^h(M) \ear\in\cR^{n_y \times {\bar \jmath}},
\eean
we get the horizontal subspace equations
\bean \label{bYph}
\bY_p^h & = & \Gamma_i^h\widehat{\bX}_p^h + \Gamma_i^{vh}\widehat{\bX}_p^{vh} + K_i^{h}\bE_p^h \\ \label{bYfh}
\bY_f^h & = & \Gamma_i^h\widehat{\bX}_f^h + \Gamma_i^{vh}\widehat{\bX}_f^{vh} + K_i^{h}\bE_f^h \\ \label{bXfh}
\widehat{\bX}_f^h & = & A_1^i\widehat{\bX}_p^h + \Phi_i^{vh}\widehat{\bX}_p^{vh} + {\cal L}_i^h\bE_p^h.
\eean
\subsection{Propagating the Vertical Hankel State Matrices}
We will now propagate the state equation \eqref{xv} backward until we reach the initial vertical states. By assuming
zero initial vertical states, then the remaining
vertical states are a function of the innovations and horizontal states only. Since \eqref{Xpv} and \eqref{Xfv} are Hankel matrices, we need to convert \eqref{xv} into a pair of past and future Hankel type matrix equations.
This is rather straightforward since \eqref{Xpv} and \eqref{Xfv} have partial horizontal dynamics $($i.e., only through $\widehat{X}_p^h(k)$ and $\widehat{X}_f^h(k))$.
Thus, by substituting $\widehat{x}_{r,s}^h$, $\widehat{x}_{r,s}^v$, and $e_{r,s}$ in \eqref{xv} by their matrix equivalents, $\{\widehat{X}_p^h(k),\widehat{X}_p^{vh}(k),E_p^h(k)\}$
and $\{\widehat{X}_f^h(k),\widehat{X}_f^{vh}(k),E_f^h(k)\}$, we obtain, respectively, the past and future vertical state equations given by
\bean \label{xs1}
\widehat{X}_p^{vh}(k+1) & = & \Theta_i^{vh}\widehat{X}_p^h(k) + A_i^{vh}\widehat{X}_p^{vh}(k) + K_i^{vh}E_p^h(k) \\ \label{xs2}
\widehat{X}_f^{vh}(k+1) & = & \Theta_i^{vh}\widehat{X}_f^h(k) + A_i^{vh}\widehat{X}_f^{vh}(k) + K_i^{vh}E_f^h(k),
\eean
where
\bean
\Theta_i^{vh} & \triangleq & \lba{c} A_3 \\ A_3A_1 \\ \vdots \\ A_3A_1^{i-1} \ear \;=\; \left(I_i \otimes A_3\right) \cdot \Theta_i^h
\in\cR^{n_vi \times n_h} \\
A_i^{vh} & \triangleq & \lba{cccc} A_4 & & & \\ A_3A_2 & A_4 & & \\
\vdots & \vdots & \ddots & \\
A_3A_1^{i-2}A_2 & A_3A_1^{i-3}A_2 & \cdots & A_4 \ear \in\cR^{n_vi \times n_vi} \nonumber \\
& = & (I_i \otimes A_3)G_{A_1}(I_i \otimes A_2) + (I_i \otimes A_4) \\
K_i^{vh} & \triangleq & \lba{cccc} K_2 & & & \\ A_3K_1 & K_2 & \\
\vdots & \vdots & \ddots & \\
A_3A_1^{i-2}K_1 & A_3A_1^{i-3}K_1 & \cdots & K_2 \ear \in\cR^{n_vi \times n_y i} \nonumber \\
& = & (I_i \otimes A_3)G_{A_1}(I_i \otimes K_1) + (I_i \otimes K_2).
\eean
Let us now solve \eqref{xs1} and \eqref{xs2} recursively for $k=0,1,\ldots,M$ as follows:
\bean
\hspace{-10mm} \widehat{X}_p^{vh}(0) & = & \widehat{X}_p^{vh}(0) \nonumber \\
\hspace{-10mm} \widehat{X}_f^{vh}(0) & = & \widehat{X}_f^{vh}(0) \nonumber \\
\hspace{-10mm} \widehat{X}_p^{vh}(1) & = & \Theta_i^{vh}\widehat{X}_p^h(0) + A_i^{vh}\widehat{X}_p^{vh}(0) + K_i^{vh}E_p^h(0) \nonumber \\
\hspace{-10mm} \widehat{X}_f^{vh}(1) & = & \Theta_i^{vh}\widehat{X}_f^h(0) + A_i^{vh}\widehat{X}_f^{vh}(0) + K_i^{vh}E_f^h(0) \nonumber \\
\hspace{-10mm} \widehat{X}_p^{vh}(2) & = & \Theta_i^{vh}\widehat{X}_p^h(1) + A_i^{vh}\Theta_i^{vh}\widehat{X}_p^h(0)
+ (A_i^{vh})^2\widehat{X}_p^{vh}(0) +
K_i^{vh}E_p^h(1) + A_i^{vh}K_i^{vh}E_p^h(0) \nonumber \\
\hspace{-10mm} \widehat{X}_f^{vh}(2) & = & \Theta_i^{vh}\widehat{X}_f^h(1) + A_i^{vh}\Theta_i^{vh}\widehat{X}_f^h(0)
+ (A_i^{vh})^2\widehat{X}_f^{vh}(0) +
K_i^{vh}E_f^h(1) + A_i^{vh}K_i^{vh}E_f^h(0)  \nonumber \\
& \vdots & \nonumber \\
\hspace{-10mm} \widehat{X}_p^{vh}(M) & = & (A_i^{vh})^{M}\widehat{X}_p^{vh}(0) + \sum_{k=0}^{M-1}(A_i^{vh})^{M-k-1}\Theta_i^{vh}\widehat{X}_p^h(k) + \sum_{k=0}^{M-1}(A_i^{vh})^{M-k-1}K_i^{vh}E_p^h(k) \\
\hspace{-10mm} \widehat{X}_f^{vh}(M) & = & (A_i^{vh})^{M}\widehat{X}_f^{vh}(0) + \sum_{k=0}^{M-1}(A_i^{vh})^{M-k-1}\Theta_i^{vh}\widehat{X}_f^h(k) + \sum_{k=0}^{M-1}(A_i^{vh})^{M-k-1}K_i^{vh}E_f^h(k).
\eean
\noindent Now we use $\{E^h_p(k),E^h_f(k),\widehat{X}^h_p(k),\widehat{X}^h_f(k)\}$ for $k=0,1,\ldots,M$ to construct
upper triangular block Toeplitz matrices such as
\bean \label{Ehp}
\bE^{\star}_{p} & \triangleq & \lba{cccc} E^h_{p}(0) & E^h_{p}(1) & \cdots & E^h_{p}(M) \\
& E^h_{p}(0) & \cdots & E^h_{p}(M-1) \\
& & \ddots & \vdots \\ & & & E^h_{p}(0) \ear\in\cR^{\bar{n}_yi \times {\bar \jmath}} \\ \label{Ehf}
\bE^{\star}_{f} & \triangleq & \lba{cccc} E^h_{f}(0) & E^h_{f}(1) & \cdots & E^h_{f}(M) \\
& E^h_{f}(0) & \cdots & E^h_{f}(M-1) \\
& & \ddots & \vdots \\ & & & E^h_{f}(0) \ear\in\cR^{\bar{n}_yi \times {\bar \jmath}} \\ \label{Xhp}
\widehat{\bX}^{\star}_{p} & \triangleq & \lba{cccc} \widehat{X}^h_{p}(0) & \widehat{X}^h_{p}(1) & \cdots & \widehat{X}^h_{p}(M) \\
& \widehat{X}^h_{p}(0) & \cdots & \widehat{X}^h_{p}(M-1) \\
& & \ddots & \vdots \\ & & & \widehat{X}^h_{p}(0) \ear\in\cR^{\bar{n}_h \times {\bar \jmath}} \\ \label{Xhf}
\widehat{\bX}^{\star}_{f} & \triangleq & \lba{cccc} \widehat{X}^h_{f}(0) & \widehat{X}^h_{f}(1) & \cdots & \widehat{X}^h_{f}(M) \\
& \widehat{X}^h_{f}(0) & \cdots & \widehat{X}^h_{f}(M-1) \\
& & \ddots & \vdots \\ & & & \widehat{X}^h_{f}(0) \ear\in\cR^{\bar{n}_h \times {\bar \jmath}}.
\eean
Notice that \eqref{Ehp} -- \eqref{Ehf} contain block Hankel entries, thus are
block Toeplitz with Hankel blocks $($BTHB$)$. Finally, we define the controllability-like matrices
\bean
\bA_M^{vh} & \triangleq & \lba{c|c|c|c} \Theta_i^{vh} & A_i^{vh}\Theta_i^{vh} & \cdots & (A_i^{vh})^{M-1}\Theta_i^{vh} \ear\in\cR^{n_vi \times n_hM} \\
\bK_M^{vh} & \triangleq & \lba{c|c|c|c} K_i^{vh} & A_i^{vh}K_i^{vh} & \cdots & (A_i^{vh})^{M-1}K_i^{vh} \ear\in\cR^{n_vi \times n_yiM}
\eean
\bean
\Delta_{\widehat{X}_{p}^{vh}(0)} & \triangleq & \lba{c|c|c|c} \widehat{X}_{p}^{vh}(0) & A_i^{vh}\widehat{X}_{p}^{vh}(0) & \cdots & (A_i^{vh})^M\widehat{X}_{p}^{vh}(0) \ear\in\cR^{n_vi \times {\bar \jmath}} \\
\Delta_{\widehat{X}_{f}^{vh}(0)} & \triangleq & \lba{c|c|c|c} \widehat{X}_{f}^{vh}(0) & A_i^{vh}\widehat{X}_{f}^{vh}(0) & \cdots & (A_i^{vh})^M\widehat{X}_{f}^{vh}(0) \ear\in\cR^{n_vi \times {\bar \jmath}}.
\eean
It can now be easily shown that the vertical states satisfy a pair of Hankel matrix equations such as
\bean \label{XXpv}
\widehat{\bX}_{p}^{vh} & = & \Delta_{\widehat{X}_{p}^{vh}(0)} + \lba{c|c} 0_{n_vi \times n_h} & \bA_M^{vh} \ear \widehat{\bX}_{p}^{\star} +
\lba{c|c} 0_{n_vi \times n_y i} & \bK_M^{vh} \ear \bE_{p}^{\star} \\ \label{XXfv}
\widehat{\bX}_{f}^{vh} & = & \Delta_{\widehat{X}_{f}^{vh}(0)} + \lba{c|c} 0_{n_vi \times n_h} & \bA_M^{vh} \ear \widehat{\bX}_{f}^{\star} +
\lba{c|c} 0_{n_vi \times n_y i} & \bK_M^{vh} \ear \bE_{f}^{\star}.
\eean
If we now assume that $\widehat{X}_{p}^{vh}(0)=\widehat{X}_{f}^{vh}(0)=0_{n_vi \times j}$, then $\Delta_{\widehat{X}_{p}^{vh}(0)}=0_{n_vi \times {\bar \jmath}}$ and $\Delta_{\widehat{X}_{f}^{vh}(0)}=0_{n_vi \times {\bar \jmath}}$. We then obtain the final expressions for $\widehat{\bX}_{p}^{vh}$ and $\widehat{\bX}_{f}^{vh}$ as
\bean \label{axpv}
\widehat{\bX}_{p}^{vh} & = & \lba{c|c} 0_{n_vi \times n_h} & \bA_M^{vh} \ear \widehat{\bX}_{p}^{\star} +
\lba{c|c} 0_{n_vi \times n_y i} & \bK_M^{vh} \ear \bE_{p}^{\star} \\ \label{axfv}
\widehat{\bX}_{f}^{vh} & = & \lba{c|c} 0_{n_vi \times n_h} & \bA_M^{vh} \ear \widehat{\bX}_{f}^{\star} +
\lba{c|c} 0_{n_vi \times n_y i} & \bK_M^{vh} \ear \bE_{f}^{\star}.
\eean
\subsection{Computing the Orthogonal Projection $\bY^h_f/\bY^h_p$}
Since the vertical state matrices are now functions of horizontal states and innovations, we can
substitute these in the horizontal state equation \eqref{bXfh}, to get
\bean \label{axfh}
\widehat{\bX}^{\star}_f & = & \bA_M^h \widehat{\bX}^{\star}_p + \bK_M^h \bE^{\star}_p,
\eean
where
\bean
\bA^h_M & \triangleq & \lba{cccccc} A_1^i & \Phi_i^{vh}\Theta_i^{vh} & \Phi_i^{vh}A_i^{vh}\Theta_i^{vh} &
\cdots & \Phi_i^{vh}(A_i^{vh})^{M-1}\Theta_i^{vh} \\
& A_1^i & \Phi_i^{vh}\Theta_i^{vh} & \cdots & \Phi_i^{vh}(A_i^{vh})^{M-2}\Theta_i^{vh} \\ & & A_1^i &
\cdots & \Phi_i^{vh}(A_i^{vh})^{M-3}\Theta_i^{vh} \\
& & & \ddots & \vdots \\
& & & & A_1^i \ear\in\cR^{\bar{n}_h \times \bar{n}_h} \\ \label{Hi_h}
\bK^h_M & \triangleq & \lba{cccccc} {\cal L}_i^h & \Phi_i^{vh}K_i^{vh} & \Phi_i^{vh}A_i^{vh}K_i^{vh} & \cdots & \Phi_i^{vh}(A_i^{vh})^{M-1}K_i^{vh} \\
& {\cal L}_i^h & \Phi_i^{vh}K_i^{vh} & \cdots & \Phi_i^{vh}(A_i^{vh})^{M-2}K_i^{vh} \\ & & {\cal L}_i^h & \cdots & \Phi_i^{vh}(A_i^{vh})^{M-3}K_i^{vh} \\
& & & \ddots & \vdots \\
& & & & {\cal L}_i^h \ear\in\cR^{\bar{n}_h \times \bar{n}_yi}.
\eean
Let us now re-visit \eqref{axpv} -- \eqref{axfv} and further substitute \eqref{axfh} in \eqref{axfv}, i.e.,
\bean
\widehat{\bX}_{p}^{vh} & = & \bQ_1 \widehat{\bX}_{p}^{\star} + \bQ_2 \bE_{p}^{\star} \\
\widehat{\bX}_{f}^{vh} & = & \bQ_1 \widehat{\bX}_{f}^{\star} + \bQ_2 \bE_{f}^{\star} \nonumber \\
& = & \bQ_1\left(\bA_M^h \widehat{\bX}^{\star}_p + \bK_M^h \bE^{\star}_p\right) + \bQ_2 \bE_{f}^{\star} \nonumber \\
& = & \bP_1\widehat{\bX}^{\star}_p + \bP_2\bE^{\star}_p + \bQ_2\bE_{f}^{\star},
\eean
where
\bean
\bQ_1 & = & \lba{c|c} 0_{n_vi \times n_h} & \bA_M^{vh} \ear\in\cR^{n_vi \times n_hM} \\
\bQ_2 & = & \lba{c|c} 0_{n_vi \times n_y i} & \bK_M^{vh}\ear\in\cR^{n_vi \times n_yiM} \\
\bP_1 & = & \lba{c|c} 0_{n_vi \times n_h} & \bA_M^{vh}\ear \cdot \bA_M^h\in\cR^{n_vi \times \bar{n}_h} \\
\bP_2 & = & \lba{c|c} 0_{n_vi \times n_h} & \bA_M^{vh}\ear \cdot \bK_M^h\in\cR^{n_vi \times \bar{n}_yi}.
\eean
We now compute the orthogonal projection $\bY^h_f/\bY^h_p$ as
\bean
\bY^h_f/\bY^h_p & = & \Gamma_i^h \widehat{\bX}^h_f/\bY^h_p + \Gamma_i^{vh}\widehat{\bX}_f^{vh}/\bY_p^h + K_i^h\bE^h_f/\bY^h_p \nonumber \\
& = & \frac{1}{\bar{\jmath}}\left(\Gamma_i^h \widehat{\bX}^h_f\left(\bY^h_p\right)^{\top} + \Gamma_i^{vh}\widehat{\bX}_f^{vh}\left(\bY_p^h\right)^{\top}
+ K_i^h\bE^h_f\left(\bY^h_p\right)^{\top}\right) \left(\bR_{pp}^{h}\right)^{-1}\bY_p^h,
\eean
where $\bR_{pp}^{h}=\frac{1}{{\bar \jmath}}\bY^h_p\left(\bY^h_p\right)^{\top}$. Furthermore, we now substitute $\left(\bY_p^h\right)^{\top}$ to get
\bea
\bY^h_f/\bY^h_p & = & \frac{1}{\bar{\jmath}}\Gamma_i^h \widehat{\bX}^h_f\left[\left(\widehat{\bX}^h_p\right)^{\top}\left(\Gamma_i^h\right)^{\top}
+ \left(\widehat{\bX}_p^{vh}\right)^{\top}\left(\Gamma_i^{vh}\right)^{\top} + \left(\bE_p^h\right)^{\top}\left(K_i^h\right)^{\top}\right]
\left(\bR_{pp}^{h}\right)^{-1}\bY_p^h \\
& & + \frac{1}{\bar{\jmath}}\Gamma_i^{vh}\widehat{\bX}_f^{vh}\left[\left(\widehat{\bX}^h_p\right)^{\top}\left(\Gamma_i^h\right)^{\top}
+ \left(\widehat{\bX}_p^{vh}\right)^{\top}\left(\Gamma_i^{vh}\right)^{\top} + \left(\bE_p^h\right)^{\top}\left(K_i^h\right)^{\top}\right]
\left(\bR_{pp}^{h}\right)^{-1}\bY_p^h \\
& & + \frac{1}{\bar{\jmath}}K_i^h\bE_f^h\left[\left(\widehat{\bX}^h_p\right)^{\top}\left(\Gamma_i^h\right)^{\top}
+ \left(\widehat{\bX}_p^{vh}\right)^{\top}\left(\Gamma_i^{vh}\right)^{\top} + \left(\bE_p^h\right)^{\top}\left(K_i^h\right)^{\top}\right]\left(\bR_{pp}^{h}\right)^{-1}\bY_p^h.
\eea
Let us now look at each term individually. We start with $\frac{1}{\bar{\jmath}}\Gamma_i^h \widehat{\bX}^h_f\left(\widehat{\bX}^h_p\right)^{\top}\left(\Gamma_i^h\right)^{\top}$
\bea
\frac{1}{\bar{\jmath}}\Gamma_i^h \widehat{\bX}^h_f\left(\widehat{\bX}^h_p\right)^{\top}\left(\Gamma_i^h\right)^{\top} & = &
\frac{1}{\bar{\jmath}}\Gamma_i^h\left(A_1^i\widehat{\bX}_p^h + \Phi_i^{vh}\widehat{\bX}_p^{vh} + {\cal L}_i^h\bE_p^h\right) \left(\widehat{\bX}^h_p\right)^{\top}\left(\Gamma_i^h\right)^{\top} \\
& = & \frac{1}{\bar{\jmath}}\Gamma_i^hA_1^i\widehat{\bX}^h_p\left(\widehat{\bX}^h_p\right)^{\top}\left(\Gamma_i^h\right)^{\top}
+ \frac{1}{\bar{\jmath}}\Gamma_i^h\Phi_i^{vh}\left(\bQ_1\widehat{\bX}_p^{\star}+\bQ_2\bE_p^{\star}\right)
\left(\widehat{\bX}^h_p\right)^{\top}\left(\Gamma_i^h\right)^{\top} \\
& & +\frac{1}{\bar{\jmath}} \Gamma_i^h{\cal L}_i^h\bE_p^h\left(\widehat{\bX}^h_p\right)^{\top}\left(\Gamma_i^h\right)^{\top} \\
& = & \Gamma_i^hA_1^iP_h\left(\Gamma_i^h\right)^{\top} + \Gamma_i^h\Phi_i^{vh} \lba{c|c} 0_{n_vi \times n_h} & \bA_M^{vh} \ear
\lba{c} P_h \\ \hline 0_{n_h \times n_h} \\ \hline \vdots \\ \hline 0_{n_h \times n_h} \ear\left(\Gamma_i^h\right)^{\top} + 0_{n_yi \times n_yi} \\
& = & \Gamma_i^hA_1^iP_h\left(\Gamma_i^h\right)^{\top}.
\eea
We continue with $\frac{1}{\bar{\jmath}}\Gamma_i^h \widehat{\bX}^h_f\left(\widehat{\bX}_p^{vh}\right)^{\top}\left(\Gamma_i^{vh}\right)^{\top}$,
\renewcommand*{\arraystretch}{1.5}
\bea
\frac{1}{\bar{\jmath}}\Gamma_i^h \widehat{\bX}^h_f\left(\widehat{\bX}_p^{vh}\right)^{\top}\left(\Gamma_i^{vh}\right)^{\top} & = &
\frac{1}{\bar{\jmath}}\Gamma_i^h\left(A_1^i\widehat{\bX}_p^h + \Phi_i^{vh}\widehat{\bX}_p^{vh} + {\cal L}_i^h\bE_p^h\right)\left(\widehat{\bX}_p^{vh}\right)^{\top}\left(\Gamma_i^{vh}\right)^{\top} \\
& = & \frac{1}{\bar{\jmath}}\Gamma_i^h A_1^i\widehat{\bX}^h_p\left(\left(\widehat{\bX}_p^{\star}\right)^{\top}\bQ_1^{\top}
+ \left(\bE_p^{\star}\right)^{\top}\bQ_2^{\top}\right)\left(\Gamma_i^{vh}\right)^{\top} \\
& & + \frac{1}{\bar{\jmath}}\Gamma_i^h\Phi_i^{vh}\widehat{\bX}_p^{vh}
\left(\widehat{\bX}_p^{vh}\right)^{\top}\left(\Gamma_i^{vh}\right)^{\top} \\
& & + \frac{1}{\bar{\jmath}}\Gamma_i^h{\cal L}_i^h\bE_p^h\left(\left(\widehat{\bX}_p^{\star}\right)^{\top}\bQ_1^{\top}
+ \left(\bE_p^{\star}\right)^{\top}\bQ_2^{\top}\right)\left(\Gamma_i^{vh}\right)^{\top}.
\eea
Thus we get,
\bea
\frac{1}{\bar{\jmath}}\Gamma_i^h \widehat{\bX}^h_f\left(\widehat{\bX}_p^{vh}\right)^{\top}\left(\Gamma_i^{vh}\right)^{\top} & = & \Gamma_i^h A_1^i \lba{c|c|c|c} P_h & 0_{n_h \times n_h} & \cdots & 0_{n_h \times n_h} \ear\lba{c}
0_{n_h \times n_vi} \\ \hline \left(\bA_M^{vh}\right)^{\top}\ear\left(\Gamma_i^{vh}\right)^{\top} + 0_{n_yi \times n_yi} \\
& & + \Gamma_i^h\Phi_i^{vh}\left(I_i \otimes P_v\right)\left(\Gamma_i^{vh}\right)^{\top}
+ 0_{n_yi \times n_yi} \\ & & + \Gamma_i^h{\cal L}_i^h \lba{c|c|c|c} R_e & 0_{n_y \times n_y} & \cdots & 0_{n_y \times n_y} \ear \lba{c}
0_{n_yi \times n_vi} \\ \hline \left(\bK_M^{vh}\right)^{\top}\ear\left(\Gamma_i^{vh}\right)^{\top} \\
& = & \Gamma_i^h\Phi_i^{vh}\left(I_i \otimes P_v\right)\left(\Gamma_i^{vh}\right)^{\top}.
\eea
Let us continue with the next term $\frac{1}{\bar{\jmath}}\Gamma_i^h \widehat{\bX}^h_f\left(\bE_p^h\right)^{\top}\left(K_i^h\right)^{\top}$
\renewcommand*{\arraystretch}{1.00}
\bea
\frac{1}{\bar{\jmath}}\Gamma_i^h \widehat{\bX}^h_f\left(\bE_p^h\right)^{\top}\left(K_i^h\right)^{\top} & = &
\frac{1}{\bar{\jmath}}\Gamma_i^h\left(A_1^i\widehat{\bX}_p^h + \Phi_i^{vh}\widehat{\bX}_p^{vh} + {\cal L}_i^h\bE_p^h\right)\left(\bE_p^h\right)^{\top}\left(K_i^h\right)^{\top} \\
& = & \frac{1}{\bar{\jmath}}\Gamma_i^hA_1^i\widehat{\bX}_p^h\left(\bE_p^h\right)^{\top}\left(K_i^h\right)^{\top} + \frac{1}{\bar{\jmath}}\Gamma_i^h\Phi_i^{vh}\left(\bQ_1\widehat{\bX}_p^{\star}
+ \bQ_2\bE_p^{\star}\right)\left(\bE_p^h\right)^{\top}\left(K_i^h\right)^{\top} \\
& & + \frac{1}{\bar{\jmath}}\Gamma_i^h{\cal L}_i^h\bE_p^h\left(\bE_p^h\right)^{\top}\left(K_i^h\right)^{\top} \\
& = & 0_{n_yi \times n_yi} + 0_{n_yi \times n_yi}+\Gamma_i^h\Phi_i^{vh}\lba{c|c} 0_{n_vi \times n_y i} & \bK_M^{vh}\ear\lba{c} R_e \\ \hline 0_{n_y \times n_y} \\ \hline \vdots \\ \hline 0_{n_y \times n_y} \ear \left(K_i^h\right)^{\top} \\
& & + \Gamma_i^h{\cal L}_i^h\left(I_i \otimes R_e\right)\left(K_i^h\right)^{\top} \\
& = & \Gamma_i^h{\cal L}_i^h\left(I_i \otimes R_e\right)\left(K_i^h\right)^{\top}.
\eea
We then continue with the term $\frac{1}{\bar{\jmath}}\Gamma_i^{vh}\widehat{\bX}_f^{vh}\left(\bY_p^h\right)^{\top}$, i.e.,
\bea
\frac{1}{\bar{\jmath}}\Gamma_i^{vh}\widehat{\bX}_f^{vh}\left(\bY_p^h\right)^{\top} & = & \frac{1}{\bar{\jmath}}\Gamma_i^{vh}\widehat{\bX}_f^{vh}
\left(\left(\widehat{\bX}_p^h\right)^{\top}\left(\Gamma_i^{h}\right)^{\top} +
\left(\widehat{\bX}_p^{vh}\right)^{\top}\left(\Gamma_i^{vh}\right)^{\top} + \left(\bE_p^h\right)^{\top}\left(K_i^h\right)^{\top}\right) \\
& = & \frac{1}{\bar{\jmath}}\Gamma_i^{vh}\widehat{\bX}_f^{vh}\left(\widehat{\bX}_p^h\right)^{\top}\left(\Gamma_i^{h}\right)^{\top}
+ \frac{1}{\bar{\jmath}}\Gamma_i^{vh}\widehat{\bX}_f^{vh}\left(\widehat{\bX}_p^{vh}\right)^{\top}\left(\Gamma_i^{vh}\right)^{\top}
+ \frac{1}{\bar{\jmath}}\Gamma_i^{vh}\widehat{\bX}_f^{vh}\left(\bE_p^{h}\right)^{\top}\left(K_i^{h}\right)^{\top} \\
& = & \frac{1}{\bar{\jmath}}\Gamma_i^{vh}\left(\bP_1\widehat{\bX}_p^{\star}+\bP_2\bE_p^{\star}+\bQ_2\bE_f^{\star}\right)
\left(\widehat{\bX}_p^h\right)^{\top}\left(\Gamma_i^{h}\right)^{\top}
+ \frac{1}{\bar{\jmath}}\Gamma_i^{vh}\widehat{\bX}_f^{vh}\left(\widehat{\bX}_p^{vh}\right)^{\top}\left(\Gamma_i^{vh}\right)^{\top} \\
& & + \frac{1}{\bar{\jmath}}\Gamma_i^{vh}\left(\bP_1\widehat{\bX}_p^{\star}+\bP_2\bE_p^{\star}+\bQ_2\bE_f^{\star}\right)
\left(\bE_p^{h}\right)^{\top}\left(K_i^{h}\right)^{\top} \\
& = & \Gamma_i^{vh}\lba{c|c} 0_{n_vi \times n_h} & \bA_M^{vh}\ear \cdot \bA_M^h \lba{c} P_h \\
\hline 0_{n_h \times n_h} \\ \hline \vdots \\ \hline 0_{n_h \times n_h} \ear \left(\Gamma_i^h\right)^{\top} + 0_{n_yi \times n_yi} + 0_{n_yi \times n_yi} \\
& & + \frac{1}{\bar{\jmath}}\Gamma_i^{vh}\widehat{\bX}_f^{vh}\left(\widehat{\bX}_p^{vh}\right)^{\top}\left(\Gamma_i^{vh}\right)^{\top} + 0_{n_yi \times n_yi}
\\ & & + \Gamma_i^{vh}\lba{c|c} 0_{n_vi \times n_h} & \bA_M^{vh}\ear \cdot \bK_M^h
\lba{c} R_e \\ \hline 0_{n_y \times n_y} \\ \hline \vdots \\ \hline 0_{n_y \times n_y} \ear \left(K_i^h\right)^{\top} + 0_{n_yi \times n_yi} \\
& = & \frac{1}{\bar{\jmath}}\Gamma_i^{vh}\widehat{\bX}_f^{vh}\left(\widehat{\bX}_p^{vh}\right)^{\top}\left(\Gamma_i^{vh}\right)^{\top}.
\eea
Finally, the last term $\frac{1}{\bar{\jmath}}K_i^h\bE^h_f\left(\bY^h_p\right)^{\top}$ is zero since the future innovations are
uncorrelated with the past data.

Now collecting all terms, we obtain
\bea
\bY^h_f/\bY^h_p & = & \Gamma_i^h\left(A_1^iP_h\left(\Gamma_i^h\right)^{\top}+\Phi_i^{vh}\left(I_i \otimes P_v\right)
\left(\Gamma_i^{vh}\right)^{\top} + {\cal L}_i^h\left(I_i \otimes R_e\right)\left(K_i^h\right)^{\top}\right) \left(\bR_{pp}^{h}\right)^{-1}\bY_p^h \\
&& + \frac{1}{\bar{\jmath}}\Gamma_i^{vh}\widehat{\bX}_f^{vh}\left(\widehat{\bX}_p^{vh}\right)^{\top}\left(\Gamma_i^{vh}\right)^{\top}
\left(\bR_{pp}^{h}\right)^{-1}\bY_p^h.
\eea
One can show that
\bea
\Delta_i^h & = & A_1^iP_h\left(\Gamma_i^h\right)^{\top} + \Phi_i^{vh}\left(I_{i} \otimes P_v\right)\left(\Gamma_i^{vh}\right)^{\top}
+ {\cal L}_i^h\left(I_i \otimes R_e\right)\left(K_i^h\right)^{\top} \\
& = & \lba{c|c|c|c} A_1^{i-1}G_1 & A_1^{i-1}G_1 & \cdots & G_1 \ear.
\eea
Therefore, we have
\bea
\bY^h_f/\bY^h_p & = & \Gamma_i^h\cdot \Delta_i^h \left(\bR_{pp}^{h}\right)^{-1}\bY_p^h +
\frac{1}{\bar{\jmath}}\Gamma_i^{vh}\widehat{\bX}_f^{vh}\left(\widehat{\bX}_p^{vh}\right)^{\top}\left(\Gamma_i^{vh}\right)^{\top}
\left(\bR_{pp}^{h}\right)^{-1}\bY_p^h.
\eea
We define the bias term as
\bea
\mbox{bias} & = & \frac{1}{\bar{\jmath}}\Gamma_i^{vh}\widehat{\bX}_f^{vh}\left(\widehat{\bX}_p^{vh}\right)^{\top}\left(\Gamma_i^{vh}\right)^{\top}
\left(\bR_{pp}^{h}\right)^{-1}\bY_p^h.
\eea
We now need to find a closed form expression for $\frac{1}{\bar{\jmath}}\widehat{\bX}_f^{vh}\left(\widehat{\bX}_p^{vh}\right)^{\top}$.
For this, we will use equations \eqref{Xf1} and \eqref{xs1} -- \eqref{xs2}, along the following state estimate covariance equations
\bea
P_{h} & = & A_1P_hA_1^{\top}+A_2P_vA_2^{\top}+K_1R_eK_1^{\top} \\
P_{hv} & = & A_1P_hA_3^{\top}+A_2P_vA_4^{\top}+K_1R_eK_2^{\top} \\
P_{vh} & = & A_3P_hA_1^{\top}+A_2P_vA_2^{\top}+K_2R_eK_1^{\top} \\
P_{v} & = & A_3P_hA_3^{\top}+A_4P_vA_4^{\top}+K_2R_eK_2^{\top}.
\eea
One can easily prove the following results.
\bea
\mathcal{P}_0 & = & \Theta_i^{vh}\cdot \left[A_1^iP_h\left(\Theta_i^{vh}\right)^{\top}+\Phi^{vh}\left(I_i \otimes P_v\right)\left(A_i^{vh}\right)^{\top}
+\mathcal{L}_i^{h}\left(I_i \otimes R_e\right)\left(K_i^{vh}\right)^{\top}\right] \\
& = & \Theta_i^{vh}\Phi_i^{h}\left(I_i \otimes P_{hv}\right) \\
\mathcal{Q} & = & \Theta_i^{vh}P_h\left(\Theta_i^{vh}\right)^{\top}+A_i^{vh}\left(I_i \otimes P_v\right)\left(A_i^{vh}\right)^{\top}
+K_i^{h}\left(I_i \otimes R_e\right)\left(K_i^{vh}\right)^{\top} \\
& = & \left(I_i \otimes P_{v}\right)+\underbrace{\left(I_i \otimes A_3\right)G_{A_1}\left(I_i \otimes P_{hv}\right) + \left(I_i \otimes P_{vh}\right)
G_{A_1}^{\top}\left(I_i \otimes A_3^{\top}\right)}_{\mathcal{Q}_0}.
\eea
Now, since $\frac{1}{\bar{\jmath}}\widehat{\bX}_f^{vh}\left(\widehat{\bX}_p^{vh}\right)^{\top}$ can be represented as
\bea
\frac{1}{\bar{\jmath}}\widehat{\bX}_f^{vh}\left(\widehat{\bX}_p^{vh}\right)^{\top} & = &
\frac{1}{\bar{\jmath}}\sum_{k=0}^{M}\widehat{X}_f^{vh}(k)\left(\widehat{X}_p^{vh}(k)\right)^{\top}
\eea
Let us consider each product term for $k=0,1,\ldots,M$. Starting with $k=0$, we have that
\bea
\frac{1}{\bar{\jmath}}\widehat{X}_f^{vh}(0)\left(\widehat{X}_p^{vh}(0)\right)^{\top} & = & 0_{n_vi \times n_vi}.
\eea
Continuing with $k=1$, we have
\renewcommand*{\arraystretch}{1.25}
\bea
\widehat{X}_p^{vh}(1) & = & \Theta_i^{vh}\widehat{X}_p^h(0) + A_i^{vh}\widehat{X}_p^{vh}(0) + K_i^{vh}E_p^h(0) \\
\widehat{X}_f^{vh}(1) & = & \Theta_i^{vh}\widehat{X}_f^h(0) + A_i^{vh}\widehat{X}_f^{vh}(0) + K_i^{vh}E_f^h(0) \\
& = & \Theta_i^{vh}A_1^i\widehat{X}_p^h(0)+\Theta_i^{vh}\Phi_i^{vh}\widehat{X}_p^{vh}(0)+\Theta_i^{vh}\mathcal{L}_i^hE_p^h(0)
+ A_i^{vh}\widehat{X}_f^{vh}(0) + K_i^{vh}E_f^h(0).
\eea
Now computing the covariance $\frac{1}{\bar{\jmath}}\widehat{X}_f^{vh}(1)\left(\widehat{X}_p^{vh}(1)\right)^{\top}$, we get
\bea
\frac{1}{\bar{\jmath}}\widehat{X}_f^{vh}(1)\left(\widehat{X}_p^{vh}(1)\right)^{\top} & = &
\frac{1}{\bar{\jmath}}\left(\Theta_i^{vh}\widehat{X}_f^h(0) + A_i^{vh}\widehat{X}_f^{vh}(0) + K_i^{vh}E_f^h(0)\right)
\left(\left(\widehat{X}_p^h(0)\right)^{\top}\left(\Theta_i^{vh}\right)^{\top} \right. \\ & & \left. +
\left(\widehat{X}_p^{vh}(0)\right)^{\top}\left(A_i^{vh}\right)^{\top} + \left(E_p^h(0)\right)^{\top}\left(K_i^{vh}\right)^{\top}\right) \\
& = & \Theta_i^{vh}\left[A_1^iP_h\left(\Theta_i^{vh}\right)^{\top} + \Phi_i^{vh}\left(I_i \otimes P_v\right)\left(A_i^{vh}\right)^{\top}
+ \mathcal{L}_i^h\left(I_i \otimes R_e\right)\left(K_i^{vh}\right)^{\top}\right] \\
& = & \Theta_i^{vh}\Phi_i^{h}\left(I_i \otimes P_{hv}\right) \\
& = & \mathcal{P}_0.
\eea
Likewise, for $k=2$, we have
\bea
\widehat{X}_p^{vh}(2) & = & \Theta_i^{vh}\widehat{X}_p^h(1) + A_i^{vh}\Theta_i^{vh}\widehat{X}_p^h(0)
+ (A_i^{vh})^2\widehat{X}_p^{vh}(0) + K_i^{vh}E_p^h(1) + A_i^{vh}K_i^{vh}E_p^h(0) \\
\widehat{X}_f^{vh}(2) & = & \Theta_i^{vh}\widehat{X}_f^h(1) + A_i^{vh}\Theta_i^{vh}\widehat{X}_f^h(0)
+ (A_i^{vh})^2\widehat{X}_f^{vh}(0) + K_i^{vh}E_f^h(1) + A_i^{vh}K_i^{vh}E_f^h(0) \\
& = & \Theta_i^{vh}A_1^i\widehat{X}_p^h(1) + \Theta_i^{vh}\Phi_i^{vh}\Theta_i^{vh}\widehat{X}_p^h(0)
+ \Theta_i^{vh}\Phi_i^{vh}A_i^{vh}\widehat{X}_p^{vh}(0) + \Theta_i^{vh}\Phi_i^{vh}K_i^{vh}E_p^{h}(0) \\
& & + \Theta_i^{vh}\mathcal{L}_i^{h}E_p^{h}(0) + A_i^{vh}\Theta_i^{vh}A_1^i\widehat{X}_p^h(0)
+ A_i^{vh}\Theta_i^{vh}\Phi_i^{vh}\widehat{X}_p^{vh}(0) + A_i^{vh}\Theta_i^{vh}\mathcal{L}_i^{h}E_p^{h}(0) \\
& & + \left(A_i^{vh}\right)^2\widehat{X}_f^{vh}(0)+A_i^{vh}K_i^{vh}E_f^h(0) + K_i^{vh}E_f^h(1)
\eea
Now computing the covariance $\frac{1}{\bar{\jmath}}\widehat{X}_f^{vh}(2)\left(\widehat{X}_p^{vh}(2)\right)^{\top}$, we get
\bea
\frac{1}{\bar{\jmath}}\widehat{X}_f^{vh}(2)\left(\widehat{X}_p^{vh}(2)\right)^{\top} & = &
\Theta_i^{vh}\left[A_1^iP_h\left(\Theta_i^{vh}\right)^{\top}+\Phi^{vh}\mathcal{Q}\left(A_i^{vh}\right)^{\top}
+\mathcal{L}_i^{h}\left(I_i \otimes R_e\right)\left(K_i^{vh}\right)^{\top}\right] \\
& & + A_i^{vh}\Theta_i^{vh}\Phi_i^{h}\left(I_i \otimes P_{hv}\right)\left(A_i^{vh}\right)^{\top} \\
& = & \Theta_i^{vh}\left[A_1^iP_h\left(\Theta_i^{vh}\right)^{\top}+\Phi^{vh}\left(I_i \otimes P_v\right)\left(A_i^{vh}\right)^{\top}
+\mathcal{L}_i^{h}\left(I_i \otimes R_e\right)\left(K_i^{vh}\right)^{\top}\right] \\
& & +\Theta_i^{vh}\Phi_i^{vh}\mathcal{Q}_0\left(A_i^{vh}\right)^{\top}
+ A_i^{vh}\Theta_i^{vh}\Phi_i^{h}\left(I_i \otimes P_{hv}\right)\left(A_i^{vh}\right)^{\top} \\
& = & \Theta_i^{vh}\Phi_i^{h}\left(I_i \otimes P_{hv}\right) +\Theta_i^{vh}\Phi_i^{vh}\mathcal{Q}_0\left(A_i^{vh}\right)^{\top} + A_i^{vh}\Theta_i^{vh}\Phi_i^{h}\left(I_i \otimes P_{hv}\right)\left(A_i^{vh}\right)^{\top} \\
& = & \mathcal{P}_0 + A_i^{vh}\mathcal{P}_0\left(A_i^{vh}\right)^{\top} + \Theta_i^{vh}\Phi_i^{vh}\mathcal{Q}_0\left(A_i^{vh}\right)^{\top}.
\eea
Continuing further, for $k \ge 3$, we obtain the general expression for
$\frac{1}{\bar{\jmath}}\Gamma_i^{vh}\widehat{\bX}_f^{vh}\left(\widehat{\bX}_p^{vh}\right)^{\top}\left(\Gamma_i^{vh}\right)^{\top}$
as
\bean \label{bias}
\frac{1}{\bar{\jmath}}\Gamma_i^{vh}\widehat{\bX}_f^{vh}\left(\widehat{\bX}_p^{vh}\right)^{\top}\left(\Gamma_i^{vh}\right)^{\top}
& = & \sum_{k=0}^{M-\ell}\sum_{\ell=1}^{M}(M-\ell - k+1)\Gamma_i^{vh}\left(A_i^{vh}\right)^{k}\mathcal{P}_0\left(\left(A_i^{vh}\right)^{k}\right)^{\top}
\left(\Gamma_i^{vh}\right)^{\top}\nonumber \\
& & + \sum_{k=0}^{M-\ell-1}\sum_{\ell=1}^{M-1}(M-\ell - k)\Gamma_i^{vh}\left(A_i^{vh}\right)^{\ell -1}\Theta_i^{vh}\Phi_i^{vh}\left(A_i^{vh}\right)^{k} \nonumber \\
& & \times \mathcal{Q}_0\left(\left(A_i^{vh}\right)^{\ell +k}\right)^{\top}\left(\Gamma_i^{vh}\right)^{\top}.
\eean
Analyzing \eqref{bias} one can see that each term in the first sum is a function of $\mathcal{P}_0$ and each term in the second sum is a
function of $\mathcal{Q}_0$, both of which are functions of $P_{hv}$ and/or $P_{vh}$, which by \eqref{Pmat} are zero matrices.
Thus, we conclude that the bias term is zero. Thus, 
\bea
\mbox{bias} & = & \frac{1}{\bar{\jmath}}\Gamma_i^{vh}\widehat{\bX}_f^{vh}\left(\widehat{\bX}_p^{vh}\right)^{\top}\left(\Gamma_i^{vh}\right)^{\top}\left(\bR_{pp}^{h}\right)^{-1}\bY_p^h \;=\; 0_{n_yi \times j_h}
\eea
and 
\bea
\bY^h_f/\bY^h_p & = & \Gamma_i^h\cdot \Delta_i^h \left(\bR_{pp}^{h}\right)^{-1}\bY_p^h \\
& = & \Gamma_i^h\cdot {\widehat \bX}^h_{f}.
\eea
\subsection{Improving The State Estimates}
Since the orthogonal projection is not exact, there is a small bias introduced that may affect the identification of
the system parameters. Despite the fact that the bias is rather small, one can iterate the procedure in order to improve
the state estimates and eliminate the bias. We now propose an oblique
projection approach to improve the state estimates. Along the way we also propose a procedure for computing
the initial states. We start by assuming that the vertical states\footnote{Here we assume that one can compute the entire vertical state sequence
$\widehat{x}_{r,s}^v$, for $r=0,1,\ldots,N$ and $s=0,1,\ldots,M$.} are available from
an orthogonal projection in the vertical direction, i.e., $\bY_f^v/\bY_p^v\cong\Gamma_i^v\cdot \widehat{\bX}_{f}^v$ $($see
\cite{Ramos2017a} for details$)$, where $\widehat{\bX}_{f}^v\in\mathbb{R}^{n_v\times \bar{\jmath}}$. Then we assemble the {\it vertical from
horizontal data processing} Hankel state matrix $\widehat{\bX}_{f}^{vh}\in\mathbb{R}^{n_vi\times \bar{\jmath}}$.

The second stage of the algorithm starts by
defining $\bW_p^h$ as
\renewcommand*{\arraystretch}{1.50}
\bean
\bW_p^h & = & \lba{c} \widehat{\bX}_p^{vh} \\ \hline \bY_p^h \ear\in\cR^{(n_v+n_y)i \times \bar{\jmath}}.
\eean
We then compute the RQ decomposition of the past/future data as follows:
\bean \label{RQ}
\lba{c} \widehat{\bX}_f^{vh} \\ \hline  \bW_p^h \\ \hline \bY_f^h \ear & = & \lba{c|c|c} R_{11} & & \\ \hline R_{21} & R_{22} & \\ \hline
R_{31} & R_{32} & R_{33} \ear \lba{c} Q_1^{\top} \\ \hline Q_2^{\top} \\ \hline Q_3^{\top} \ear,
\eean
where $R_{11}\in\mathbb{R}^{n_vi \times n_vi}$, $R_{21}\in\mathbb{R}^{(n_v+n_y)i \times n_vi}$,
$R_{22}\in\mathbb{R}^{(n_v+n_y)i \times (n_v+n_y)i}$, $R_{31}\in\mathbb{R}^{n_yi \times n_vi}$, $R_{32}\in\mathbb{R}^{n_yi \times (n_v+n_y)i}$,
$R_{33}\in\mathbb{R}^{n_yi \times n_yi}$, $Q_{1}\in\mathbb{R}^{j \times n_vi}$, $Q_{2}\in\mathbb{R}^{j \times (n_v+n_y)i}$, and
$Q_{3}\in\mathbb{R}^{j \times n_yi}$.

From \eqref{RQ} one can find an expression for $\bY_f^h$ using the R and Q parameters, along with using $Q_1^{\top}=R_{11}^{-1}\widehat{\bX}_f^{vh}$ and
$Q_2^{\top}=R_{22}^{-1}\left(\bW_p^h - R_{21}Q_1^{\top}\right)$. That is,
\renewcommand*{\arraystretch}{1.00}
\bea
\bY_f^h & = & R_{31}Q_1^{\top} + R_{32}Q_2^{\top} + R_{33}Q_3^{\top} \\
& = & R_{31}Q_1^{\top} + R_{32}R_{22}^{-1}\left(\bW_p^h - R_{21}Q_1^{\top}\right) + R_{33}Q_3^{\top} \\
& = & R_{32}R_{22}^{-1}\bW_p^h + \left(R_{31}-R_{32}R_{22}^{-1}R_{21}\right)R_{11}^{-1}\cdot \widehat{\bX}_f^{vh} + R_{33}Q_3^{\top} \\
& = & \Gamma_i^h \cdot \widehat{\bX}_f^h + \Gamma_i^{vh}\cdot \widehat{\bX}_f^{vh} + K_i^h\cdot \bE_f^h.
\eea
It is now clearly evident that
\bea
\Gamma_i^h \cdot \widehat{\bX}_f^h & = & R_{32}R_{22}^{-1}\bW_p^h \\
\Gamma_i^{vh}\cdot \widehat{\bX}_f^{vh} & = & \left(R_{31}-R_{32}R_{22}^{-1}R_{21}\right)R_{11}^{-1}\cdot \widehat{\bX}_f^{vh} \\
K_i^h\cdot \bE_f^h & = & R_{33}Q_3^{\top}.
\eea
Without computing the system parameters, our aim here is to compute $\Gamma_i^h$, then $\Gamma_i^{vh}$ and $K_i^h$
with the right lower triangular Toeplitz structure. Computing $\Gamma_i^h$ is straight forward, thus we assume it is
already known. We now concentrate on computing $\Gamma_i^{vh}$ and $K_i^h$.

Using the Lower Triangular Toeplitz System Solver $($LTTSS$)$ procedure in \cite{Ramos2016b}, we compute $\Gamma_i^{vh}$ by solving the linear system of equations
\bea
I_{n_yi} \cdot \Gamma_i^{vh} \cdot R_{11} & = & \left(R_{31}-R_{32}R_{22}^{-1}R_{21}\right),
\eea
subject to $\Gamma_i^{vh}$ being lower triangular Toeplitz. The solution is
\bea
\Gamma_i^{vh} & = & \mbox{LTTSS}\{I_{n_yi},R_{11},\left(R_{31}-R_{32}R_{22}^{-1}R_{21}\right),n_y,n_v,i\}.
\eea
We now define $\bE_f$ and $\bE_{f_1}$ as
\bea
\bE_f & = & R_{33}Q_3^{\top} \;=\; \lba{c|c|c|c} E_f(0) & E_f(1) & \cdots & E_f(M) \ear\in\mathbb{R}^{n_yi \times \bar{\jmath}} \\
\bE_{f_1} & = & \lba{c|c|c|c} E_{f_1}(0) & E_{f_1}(1) & \cdots & E_{f_1}(M) \ear\in\mathbb{R}^{n_yi \times (j-i+1)(M+1)},
\eea
where
\bea
E_f(k) & = & \lba{ccccc|ccc} e^k_{0,0} & e^k_{0,1} & e^k_{0,2} & \cdots & e^k_{0,j-i} & e^k_{0,j-i+1} & \cdots & e^k_{0,j-1} \\
e^k_{1,0} & e^k_{1,1} & e^k_{1,2} & \cdots & e^k_{1,j-i} & e^k_{1,j-i+1} & \cdots & e^k_{1,j-1} \\
e^k_{2,0} & e^k_{2,1} & e^k_{2,2} & \cdots & e^k_{2,j-i} & e^k_{2,j-i+1} & \cdots & e^k_{2,j-1} \\
\vdots & \vdots & \vdots & \ddots & \vdots & \ddots & \vdots & \vdots \\
e^k_{i-1,0} & e^k_{i-1,1} & e^k_{i-1,2} & \cdots & e^k_{i-1,j-i} & e^k_{i-1,j-i+1} & \cdots & e^k_{i-1,j-1} \ear \\
& = & \lba{c|c} E_{f_1}(k) & \times \ear
\eea
and $\times$ denotes a matrix that is not relevant to the discussion.
Furthermore, since the main diagonal blocks of $K_i^h$ are all equal to $I_{n_y}$ and all
elements above the main diagonal blocks are $0_{n_y\times n_y}$, we observe that the first $n_y$ rows
of $\bE_f$ contains a sequence of innovations, from which $K_i^h$ can be computed. That is,
\bea
\lba{c|c|c|c} I_{n_y} & 0_{n_y \times n_y} &
\cdots & 0_{n_y \times n_y} \ear \bE_f & = & \lba{c|c|c|c} \bbe_0(0) & \bbe_0(1) & \cdots & \bbe_0(M) \ear\in\mathbb{R}^{n_y \times \bar{\jmath}},
\eea
where
\bea
\bbe_0(k) & = & \lba{c|c|c|c} e^k_{0,0} & e^k_{0,1} & \cdots & e^k_{0,j-1} \ear\;\cong
\lba{c|c|c|c} e_{i,k} & e_{i+1,k} & \cdots & e_{i+j-1,k} \ear\in\mathbb{R}^{n_y \times j}.
\eea
Let us now form the array of Hankel matrices using $\bbe_0(k)$, for $k=0,1,\ldots,M$, i.e.,
\bea
\bE_{f_2} & = & \lba{c|c|c|c} E_{f_2}(0) & E_{f_2}(1) & \cdots & E_{f_2}(M)\ear\in\mathbb{R}^{n_yi \times (j-i+1)(M+1)},
\eea
where
\bea
E_{f_2}(k) & = & \lba{ccccc} e^k_{0,0} & e^k_{0,1} & e^k_{0,2} & \cdots & e^k_{0,j-i} \\
e^k_{0,1} & e^k_{0,2} & e^k_{0,3} & \cdots & e^k_{0,j-i+1} \\
e^k_{0,2} & e^k_{0,3} & e^k_{0,4} & \cdots & e^k_{0,j-i+2} \\
\vdots & \vdots & \vdots & \iddots & \vdots \\
e^k_{0,i-1} & e^k_{0,i} & e^k_{0,i+1} & \cdots & e^k_{0,j-1} \ear\in\mathbb{R}^{n_yi \times (j-i+1)}.
\eea
Notice that if we knew $\bE_f^h$, then $E_{f_2}(k)$ would be the first $j-i+1$ columns of $E_f^h(k)$, for $k=0,1,\ldots,M$. That is,
\bea
E_f^h(k) & = & \lba{ccccc|ccc} e_{i,k} & e_{i+1,k} & e_{i+2,k} & \cdots & e_{j,k} & e_{j+1,k} & \cdots & e_{i+j-1,k} \\
e_{i+1,k} & e_{i+2,k} & e_{i+3,k} & \cdots & e_{j+1,k} & e_{j+2,k} & \cdots & e_{i+j,k} \\
e_{i+2,k} & e_{i+3,k} & e_{i+4,k} & \cdots & e_{j+2,k} & e_{j+3,k} & \cdots & e_{i+j+1,k} \\
\vdots & \vdots & \vdots & \iddots & \vdots & \vdots & \iddots & \vdots \\
e_{2i-1,k} & e_{2i,k} & e_{2i+1,k} & \cdots & e_{i+j-1,k} & e_{i+j,k} & \cdots & e_{2i+j-2,k} \ear \\
& = & \lba{c|c} E_{f_2}(k) & \times \ear.
\eea
Let us now define the covariance matrices $\bV_1\in\mathbb{R}^{n_yi \times n_yi}$ and $\bV_2\in\mathbb{R}^{n_yi \times n_yi}$ as
\bea
\bV_1 & = & \frac{1}{(j-i+1)(M+1)}\bE_{f_1}\left(\bE_{f_2}\right)^{\top} \\
\bV_2 & = & \frac{1}{(j-i+1)(M+1)}\bE_{f_2}\left(\bE_{f_2}\right)^{\top}.
\eea
We can now find a relationship between $\bV_1$ and $\bV_2$ as follows:
\bea
\bV_1 & = & I_{n_y} \cdot K_i^h \cdot \bV_2,
\eea
and upon applying the vec operator on both sides, we obtain
\bean \label{Vn}
\vect\{\bV_1\} & = & \left(\bV_2^T \otimes I_{n_yi}\right) \cdot \vect\{K_i^h\}.
\eean
Now, since $K_i^h$ is a lower triangular Toeplitz matrix, $\vect\{K_i^h\}$ will
contain repeated elements. To remove these redundancies, we apply the identity
\bea
\vect\{K_i^h\} & = & {\cal F}_{K_i^h}\cdot k_i^h,
\eea
where ${\cal F}_{K_i^h}\in\mathbb{R}^{n_y^2i^2\times n_y^2i}$ is a permutation matrix with elements equal to $0$ and $1$ and $k_i^h$ contains all the elements of $K_i^h$, i.e.,
\bea
k_i^h & = & \lba{c} \vect\{k_0\} \\ \vect\{k_1\} \\ \vdots \\ \vect\{k_{i-1}\} \ear,
\eea
with
\bea
k_s & = & \left\{\ba{ll} I_{n_y}, & \mbox{if}\;s=0 \\
C_1A_1^{s-1}K_1, & \mbox{if}\;s\ge 1. \ea \right.
\eea
The $($LTTSS$)$ procedure performs this operation and assembles the full Toeplitz matrix. Thus, we obtain $K_i^h$ as
\bea
V & = & \mbox{LTTSS}\{I_{n_yi},\bV_2,\bV_1,n_y,n_y,i\} \\
K_i^h & = & V\left(I_i \otimes K_0^{-1}\right),
\eea
where $K_0\in\mathbb{R}^{n_y \times n_y}$ is the first $(n_y \times n_y)$ block of $V$.
Knowing $K_i^{h}$, we can now compute $\bE_f^{h}$ from $R_{33}Q_3^{\top}$.
For $k=0,1,\dots,M$, we need to solve for $E_f(k)$ using the Hankel System Solver $($HSS$)$ procedure
outlined in \cite{Ramos2016a}. That is, we solve the following linear system of equations
\bea
E_f(k) & = & K_i^h\cdot E_f^h(k)\cdot I_j, \; \mbox{for}\; k=0,1,\ldots,M,
\eea
subject to $E_f(k)$ being a Hankel matrix.
Upon applying the vec operator on both sides, we get
\bean \label{E_eq}
\vect\{E_f(k)\} & = & \left(I_j \otimes K_i^h\right)\vect\{E_f^h(k)\}.
\eean
However, since $E_f^h(k)$ is a Hankel matrix, $\vect\{E_f^h(k)\}$, will be inefficient for solving \eqref{E_eq}
because of the repeated elements. In order to compute the minimum number of elements from $\vect\{E_f^h(k)\}$, we
use the property
\bea
\vect\{E_f^h(k)\} & = & {\cal F}_{E_f^h(k)}\cdot e^h_f(k),
\eea
where ${\cal F}_{E_f^h(k)}\in\mathbb{R}^{n_yij \times n_y(i+j-1)}$ is a permutation matrix with elements equal to $0$ and $1$ and
\bea
e^h_f(k) & = & \lba{c} \vect\{e_{i,k}\} \\ \vect\{e_{i+1,k}\} \\ \vdots \\ \vect\{e_{2i+j-2,k}\} \ear.
\eea
The HSS procedure handles the removal of repeated elements and assembles the full Hankel matrix $E_F^h(k)$. Thus, we get
\bea
E_f^h(k) & = & \mbox{HSS}\{K_i^h,I_j,E_f(k),n_y,1,i,j\},\;\mbox{for}\;k = 0,1,\ldots,M.
\eea
Once we have all $M+1$ solutions, we can then assemble the full matrix $\bE_f^h$.

We now need to find $\widehat{\bX}_f^{vh}$ using $\{\bY_f^h,\widehat{\bX}_f^h,\bE_f^{h},\Gamma_i^h,\Gamma_i^{vh},K_i^{h}\}$. That is, let us define
\bea
\widehat{\bZ}^{vh}_f & = & \bY_f^h-\Gamma_i^h\widehat{\bX}_f^h-K_i^{h}\bE_f^{h} \;=\; \Gamma_i^{vh}\widehat{\bX}_f^{vh} \\
& = & \lba{c|c|c|c} \widehat{Z}^{vh}_f(0) & \widehat{Z}^{vh}_f(1) & \cdots & \widehat{Z}^{vh}_f(M) \ear\in\mathbb{R}^{n_yi \times \bar{\jmath}}.
\eea
Once again, we have $M+1$ systems of equations of the form
\bea
\widehat{Z}^{vh}_f(k) & = & \Gamma_i^{vh}\cdot \widehat{X}_f^{vh}(k) \cdot I_j,\;\mbox{for}\;k = 0,1,\ldots,M.
\eea
If we now apply the vec operator on both sides, we get
\bean
\vect\{\widehat{Z}^{vh}_f(k)\} & = & \left(I_j \otimes \Gamma_i^{vh}\right)\cdot {\cal F}_{\widehat{X}_f^{vh}(k)} \cdot \widehat{x}^{vh}_f(k),
\eean
where ${\cal F}_{\widehat{X}_f^{vh}(k)}\in\mathbb{R}^{n_vij \times n_v(i+j-1)}$ is a permutation matrix with elements equal to $0$ and $1$ and
\bea
\widehat{x}^{vh}_f(k) & = & \lba{c} \vect\{\widehat{x}^v_{i,k}\} \\ \vect\{\widehat{x}^v_{i+1,k}\} \\ \vdots \\
\vect\{\widehat{x}^v_{2i+j-2,k}\} \ear.
\eea
By applying the $($HSS$)$ procedure, we obtain
\bea
\widehat{X}_f^{vh}(k) & = & \mbox{HSS}\{\Gamma_i^{vh},I_j,\widehat{Z}_f^{vh}(k),n_v,1,i,j\},\;\mbox{for}\;k = 0,1,\ldots,M.
\eea

In horizontal data processing, as it relates to future data, we have two instances where we need to use the HSS procedure. This operation could be computationally expensive since it has to be done $(M+1)$ times. Nevertheless, the solution will give us the right structure
for $\bE_f^{h}$ and $\widehat{\bX}_f^{vh}$. Now that we have $\widehat{X}_f^{vh}(k)$ for $k=0,1,\ldots,M$, we can assemble the full matrix $\widehat{\bX}_f^{vh}$.

Let us now define the following matrices:
\renewcommand*{\arraystretch}{1.50}
\bea
\bT^h_2 & = & \lba{c|c|c} -\left(\Gamma_i^h\right)^{\dagger}\Gamma_i^{vh} & -\left(\Gamma_i^h\right)^{\dagger}K_i^h &
\left(\Gamma_i^h\right)^{\dagger} \ear\in\mathbb{R}^{n_h \times (n_v+2n_y)i} \\
\bH_f^h & = & \lba{c} \widehat{\bX}_f^{vh} \\ \hline \bE_f^{h} \\ \hline \bY_f^{h} \ear\in\mathbb{R}^{(n_v+2n_y)i \times \bar{\jmath}},
\eea
where $\left(\Gamma_i^h\right)^{\dagger}$ is the pseudo-inverse of $\Gamma_i^h$.
Notice that the future horizontal state matrix $\widehat{\bX}_f^{h}$ is related to $\bT_2^h$ and $\bH_f^h$ via
\bea
\widehat{\bX}_f^{h} & = & \bT_2^h\bH_f^h.
\eea

The next step is to recover $\widehat{\bX}_p^{h}$. Let us re-visit the QR decomposition of the data, i.e.,
\renewcommand*{\arraystretch}{1.25}
\bean \label{RQ}
\lba{c} \widehat{\bX}_f^{vh} \\ \hline  \widehat{\bX}_p^{vh} \\ \bY_p^h \\ \hline \bY_f^h \ear & = & \lba{c|cc|c} R_{11} & & &
\\ \hline R^1_{21} & R^1_{22} & & \\ R^2_{21} & R^2_{22} & R^3_{22} & \\ \hline
R_{31} & R^1_{32} & R^2_{32} & R_{33} \ear \lba{c} Q_1^{\top} \\ \hline Q_{21}^{\top} \\ Q_{22}^{\top} \\ \hline Q_3^{\top} \ear,
\eean
where $R_{21}^1\in\mathbb{R}^{n_vi \times n_vi}$, $R_{21}^2\in\mathbb{R}^{n_yi \times n_vi}$, $R_{22}^1\in\mathbb{R}^{n_vi \times n_vi}$,
$R_{22}^2\in\mathbb{R}^{n_yi \times n_vi}$, $R_{22}^3\in\mathbb{R}^{n_yi \times n_yi}$, $R_{32}^1\in\mathbb{R}^{n_yi \times n_vi}$,
$R_{32}^2\in\mathbb{R}^{n_yi \times n_yi}$, $Q_{21}\in\mathbb{R}^{j \times n_vi}$, and $Q_{22}\in\mathbb{R}^{j \times n_yi}$.
Then $Q_{21}^{\top}$ and $\bY_p^h$ can be expressed as
\bea
Q_{21}^{\top} & = & \left(R^1_{22}\right)^{-1}\widehat{\bX}_p^{vh} - \left(R^1_{22}\right)^{-1}R_{21}^1Q_{1}^{\top} \\
\bY_p^h & = & R_{21}^2Q_1^{\top} + R_{22}^2Q_{21}^{\top} + R_{22}^3Q_{22}^{\top} \\
& = & R_{21}^2Q_1^{\top} + R_{22}^2\left(\left(R^1_{22}\right)^{-1}\widehat{\bX}_p^{vh} - \left(R^1_{22}\right)^{-1}R_{21}^1Q_{1}^{\top}\right)
+ R_{22}^3Q_{22}^{\top} \\
& = & \left(R_{21}^2-R_{22}^2\left(R^1_{22}\right)^{-1}R_{21}^1\right)Q_{1}^{\top} + R_{22}^2\left(R^1_{22}\right)^{-1}\widehat{\bX}_p^{vh}
+ R_{22}^3Q_{22}^{\top}.
\eea
We can now isolate $R_{22}^2\left(R^1_{22}\right)^{-1}\widehat{\bX}_p^{vh}$ from the rest. That is, we define
\bea
\widehat{\bZ}_p^{vh} & = & \bY_p^h  - \left(R_{21}^2-R_{22}^2\left(R^1_{22}\right)^{-1}R_{21}^1\right)Q_{1}^{\top} - R_{22}^3Q_{22}^{\top} \;=\;
\underbrace{{R_{22}^2}\left(R^1_{22}\right)^{-1}}_{\Gamma_i^{vh}} \cdot \widehat{\bX}_p^{vh} \\
& = & \lba{c|c|c|c} \widehat{Z}_p^{vh}(0) & \widehat{Z}_p^{vh}(1) & \cdots & \widehat{Z}_p^{vh}(M) \ear
\;=\;\Gamma_i^{vh} \cdot \lba{c|c|c|c} \widehat{X}_p^{vh}(0) & \widehat{X}_p^{vh}(1) & \cdots & \widehat{X}_p^{vh}(M) \ear,
\eea
where $\widehat{\bZ}_p^{vh}\in\mathbb{R}^{n_yi \times \bar{\jmath}}$ and $\widehat{Z}_p^{vh}(k)\in\mathbb{R}^{n_yi \times j}$. It is clear that we have $M+1$ equations of the form
\bean \label{Zpp}
\widehat{Z}_p^{vh}(k) & = & \Gamma_i^{vh} \cdot \widehat{X}_p^{vh}(k) \cdot I_j,\;\mbox{for}\;k=0,1,\ldots,M.
\eean
If we apply the vec operator on both sides of \eqref{Zpp}, we get
\bean \label{Vq}
\vect\{\widehat{Z}_p^{vh}(k)\} & = & \left(I_j \otimes \Gamma_i^{vh}\right) {\cal F}_{\widehat{X}_p^{vh}(k)}\cdot \widehat{x}_p^{vh},
\eean
where ${\cal F}_{\widehat{X}_p^{vh}(k)}\in\mathbb{R}^{n_vij\times n_v(i+j-1)}$ is a permutation matrix with elements equal to $0$ and $1$
and
\bea
\widehat{x}_p^{vh}(k) & = & \lba{c} \widehat{x}^v_{0,k} \\ \widehat{x}^v_{1,k} \\ \vdots \\ \widehat{x}^v_{i+j-2,k} \ear.
\eea
Now applying the HSS procedure, which solves \eqref{Vq} and reconstructs the Hankel matrix, we can find
the individual solutions from
\bea
\widehat{X}_p^{vh}(k) & = & \mbox{HSS}\{\Gamma_i^{vh},I_j,\widehat{Z}_p^{vh}(k),n_v,1,i,j\},\;\mbox{for}\;k = 0,1,\ldots,M.
\eea
Knowing all $\widehat{X}_p^{vh}(k)$, for $k=0,1,\ldots,M$, we can now assemble the full $\widehat{\bX}_p^{vh}$ matrix.

We now let
\renewcommand*{\arraystretch}{1.00}
\bea
\bY_p & = & \left(\Gamma_i^h\right)^{\perp}\left(\bY_p^h-\Gamma_i^{vh}\cdot \widehat{\bX}_p^{vh}\right) \; = \; \underbrace{\left(\Gamma_i^h\right)^{\perp}\Gamma_i^h\widehat{\bX}_p^{h}}_{0_{(n_yi-n_h)\times \bar{\jmath}}} + \left(\Gamma_i^h\right)^{\perp}K_i^h\bE_p^{h} \\
& = & \lba{c|c|c|c} Y_p(0) & Y_p(1) & \cdots & Y_p(M) \ear
\;=\;\left(\Gamma_i^h\right)^{\perp}K_i^h\lba{c|c|c|c} E^h_p(0) & E^h_p(1) & \cdots & E^h_p(M) \ear,
\eea
where $Y_p(k)\in\mathbb{R}^{(n_yi-n_h) \times j}$ and $\left(\Gamma_i^h\right)^{\perp}\in\mathbb{R}^{(n_yi-n_h) \times n_yi}$ is the orthogonal complement of
$\Gamma_i^h$. We now have $M+1$ equations of the form
\bean \label{Ypp2}
Y_p(k) & = & \left(\Gamma_i^h\right)^{\perp}K_i^h\cdot E_p^{h}(k)\cdot I_j,
\;\mbox{for}\;k=0,1,\ldots,M.
\eean
Let us now apply the vec operator on both sides of \eqref{Ypp2}, to get
\bea
\vect\{Y_p(k)\} & = & \left(I_j \otimes \left(\Gamma_i^h\right)^{\perp}K_i^h \right)\cdot {\cal F}_{E_p^{h}(k)}\cdot e_p^h(k), \;\mbox{for}\;k=0,1,\ldots,M,
\eea
where ${\cal F}_{E_p^h(k)}\in\mathbb{R}^{n_yij \times n_y(i+j-1)}$ is a permutation matrix with elements equal to $0$ and $1$ and
\bea
e^h_p(k) & = & \lba{c} \vect\{e_{0,k}\} \\ \vect\{e_{1,k}\} \\ \vdots \\ \vect\{e_{i+j-2,k}\} \ear.
\eea
Once again, the solution can be found by applying the HSS procedure,
\bea
E^h_p(k) & = & \mbox{HSS}\{\left(\Gamma_i^h\right)^{\perp}K_i^h,I_j,Y_p(k),n_y,1,i,j\},\;\mbox{for}\;k = 0,1,\ldots,M,
\eea
from which one can assemble the full $\bE^h_p$ matrix.

Knowing $\widehat{\bX}^{vh}_p$ and $\bE^h_p$, one can now compute $\widehat{\bX}_p^{h}$ from
\renewcommand*{\arraystretch}{1.50}
\bea
\widehat{\bX}_p^{h} & = & \left(\Gamma_i^h\right)^{\dagger}\left(\bY_p^h - \Gamma_i^{vh}\widehat{\bX}_p^{vh} - K_i^h\bE_p^{h}\right) \\
& = & \lba{c|c|c} -\left(\Gamma_i^h\right)^{\dagger}\Gamma_i^{vh} & -\left(\Gamma_i^h\right)^{\dagger}K_i^h & \left(\Gamma_i^h\right)^{\dagger} \ear \lba{c} \widehat{\bX}_p^{vh} \\ \hline \bE_p^{h} \\ \hline \bY_p^{h} \ear \\
& = & \bT_2^h \bH_p^h,
\eea
where
\renewcommand*{\arraystretch}{1.50}
\bea
\bH_p^h & = & \lba{c} \widehat{\bX}_p^{vh} \\ \hline \bE_p^{h} \\ \hline \bY_p^{h} \ear\in\mathbb{R}^{(n_v+2n_y)i \times \bar{\jmath}}.
\eea
Let us now define $\bJ\in\mathbb{R}^{n_h \times (n_h+(n_v+n_y)i)}$ and $\bH\in\mathbb{R}^{(n_h+(n_v+n_y)i) \times \bar{\jmath}}$ as
\bea
\bJ & = & \lba{c|c|c} A_1^i & \Phi_i^{vh} & {\cal L}_i^h \ear
\;\;\mbox{and}\;\;\bH\;=\;\lba{c} \widehat{\bX}_p^{h} \\ \hline \widehat{\bX}_p^{vh} \\ \hline \bE_p^{h} \ear.
\eea
Then since $\widehat{\bX}_f^h=\bJ\bH$, solving for $\bJ$
would require solving a large system  of equations. In order to avoid this, we compute the
covariances $\bZ_1$ and $\bZ_2$ as
\bea
\bZ_1 & = & \frac{1}{\bar{\jmath}}\bH\bH^{\top}\in\mathbb{R}^{(n_h+(n_v+n_y)i) \times (n_h+(n_v+n_y)i)} \\
\bZ_2 & = & \frac{1}{\bar{\jmath}}\widehat{\bX}_f^h\bH^{\top}\in\mathbb{R}^{n_h \times (n_h+(n_v+n_y)i)}.
\eea
Then, the solution for $\bJ$ becomes
\bea
\bJ & = & \bZ_2\bZ^{-1}_1,
\eea
from which $A_1^i$, $\Phi_i^{vh}$, and ${\cal L}_i^h$ can be computed, i.e.,
\bea
A_1^i & = & \bJ(:,1:n_h) \\
\Phi_i^{vh} & = & \bJ(:,n_h +1:n_h+n_vi) \\
{\cal L}_i^h & = & \bJ(:,n_h+n_vi+1:n_h+(n_v+n_y)i).
\eea

Now substituting $\widehat{\bX}_p^h$ into $\widehat{\bX}_f^h$ we get
\bea
\widehat{\bX}_f^h & = & A_1^i\widehat{\bX}_p^h + \Phi_i^{vh}\widehat{\bX}_p^{vh} + {\cal L}_i^h\bE_p^h \\
& = & A_1^i\lba{c|c|c} -\left(\Gamma_i^h\right)^{\dagger}\Gamma_i^{vh} & -\left(\Gamma_i^h\right)^{\dagger}K_i^h & \left(\Gamma_i^h\right)^{\dagger} \ear \lba{c} \widehat{\bX}_p^{vh} \\ \hline \bE_p^{h} \\ \hline \bY_p^{h} \ear + \Phi_i^{vh}\widehat{\bX}_p^{vh} + {\cal L}_i^h\bE_p^h \\
& = & \lba{c|c|c} \Phi_i^{vh} - A_1^i\left(\Gamma_i^h\right)^{\dagger}\Gamma_i^{vh} & {\cal L}_i^h -A_1^i\left(\Gamma_i^h\right)^{\dagger}K_i^h &
A_1^i\left(\Gamma_i^h\right)^{\dagger} \ear \lba{c} \widehat{\bX}_p^{vh} \\ \hline \bE_p^{h} \\ \hline \bY_p^{h} \ear \\
& = & \bT_1^h \bH_p^h,
\eea
where $\bT_1^h\in\mathbb{R}^{n_h \times (n_v+2n_y)i}$ is defined as
\bea
\bT_1^h & = & \lba{c|c|c} \Phi_i^{vh} - A_1^i\left(\Gamma_i^h\right)^{\dagger}\Gamma_i^{vh} & {\cal L}_i^h -A_1^i\left(\Gamma_i^h\right)^{\dagger}K_i^h &
A_1^i\left(\Gamma_i^h\right)^{\dagger} \ear.
\eea
The final task involves computing $\widehat{\bX}_{f+}^h$ from
\bea
\widehat{\bX}_{f+}^h & = & \bT_1^h\bH_f^h \;=\;
\lba{c|c|c|c} \widehat{X}_{f+}^h(0) & \widehat{X}_{f+}^h(1) & \cdots & \widehat{X}_{f+}^h(M) \ear\in\mathbb{R}^{n_h \times \bar{\jmath}},
\eea
where
\bea
\widehat{X}_{f+}^h(k) & = & \lba{c|c|c|c} \widehat{x}_{2i,k}^h & \widehat{x}_{2i+1,k}^h & \cdots & \widehat{x}_{2i+j-1,k}^h
\ear\in\mathbb{R}^{n_h \times j}.
\eea
Now, between $\widehat{X}_{p}^h(k)$ and $\widehat{X}_{f+}^h(k)$ there is an overlap of $\{\widehat{x}_{2i,k}^h,\widehat{x}_{2i+1,k}^h,
\ldots,\widehat{x}_{j-1,k}^h\}$. That is,
\bea
\lba{c|c|c|c|c} \underbrace{\ba{c|c|c|c|c|c|c|c|c} \widehat{x}_{0,k}^h & \widehat{x}_{1,k}^h & \widehat{x}_{2,k}^h & \cdots &
\widehat{x}_{2i-1,k}^h & \widehat{x}_{2i,k}^h & \widehat{x}_{2i+1,k}^h & \cdots &
\widehat{x}_{j-1,k}^h\ea}_{\widehat{X}_{p}^h(k)} & \widehat{x}_{j,k}^h & \widehat{x}_{j+1,k}^h & \cdots & \widehat{x}_{2i+j-1,k}^h \ear \\
\eea
and
\bea
\lba{c|c|c|c|c|c} \widehat{x}_{0,k}^h & \widehat{x}_{1,k}^h & \widehat{x}_{2,k}^h & \cdots & \widehat{x}_{2i-1,k}^h
& \underbrace{\ba{c|c|c|c|c|c|c|c} \widehat{x}_{2i,k}^h & \widehat{x}_{2i+1,k}^h & \cdots &
\widehat{x}_{j-1,k}^h & \widehat{x}_{j,k}^h & \widehat{x}_{j+1,k}^h & \cdots & \widehat{x}_{2i+j-1,k}^h\ea}_{\widehat{X}_{f+}^h(k)} \ear.
\eea
Therefore,
the entire horizontal state sequence can be recovered from
\bea
\lba{c|c} \widehat{X}_p^h(:,0:2i-1)(k) & \widehat{X}_{f+}^h(k) \ear,\; \mbox{for}\; k=0,1,\dots,M.
\eea
Now form the $(n_h(N+1) \times (M+1))$ matrix of horizontal states $\widehat{\bX}^h$ by vectorizing the individual terms in
$\widehat{X}_p^h(k)$ and $\widehat{X}_{f+}^h(k)$, for $k=0,1,\ldots,M$, then recalling that $N=2i+j-2$, i.e.,
\bean
\widehat{\cal X}_p^h & \triangleq & \lba{c|c|c|c} \vect\{\widehat{X}_p^{h}(0)\} & \vect\{\widehat{X}_p^{h}(1)\} &
\cdots & \vect\{\widehat{X}_p^{h}(M)\} \ear\in\mathbb{R}^{n_hj \times (M+1)} \\
\widehat{\cal X}_{f+}^h & \triangleq & \lba{c|c|c|c} \vect\{\widehat{X}_{f+}^{h}(0)\} & \vect\{\widehat{X}_{f+}^{h}(1)\} & \cdots & \vect\{\widehat{X}_{f+}^{h}(M)\} \ear\in\mathbb{R}^{n_hj \times (M+1)},
\eean
and finally compute $\widehat{\bX}^h$ from
\bean
\widehat{\bX}^h & \triangleq & \lba{c} \widehat{\cal X}_p^h(0:2i-1,:) \\ \hline \widehat{\cal X}_{f+}^h(0:j-2,:) \ear \nonumber \\
& = & \lba{cccc} \widehat{x}_{0,0}^h & \widehat{x}_{0,1}^h & \cdots & \widehat{x}_{0,M}^h \\
\widehat{x}_{1,0}^h & \widehat{x}_{1,1}^h & \cdots & \widehat{x}_{1,M}^h \\
\vdots & \vdots & \ddots & \vdots \\
\widehat{x}_{2i-1,0}^h & \widehat{x}_{2i-1,1}^h & \cdots & \widehat{x}_{2i-1,M}^h \\ \hline
\widehat{x}_{2i,0}^h & \widehat{x}_{2i,1}^h & \cdots & \widehat{x}_{2i,M}^h \\
\widehat{x}_{2i+1,0}^h & \widehat{x}_{2i+1,1}^h & \cdots & \widehat{x}_{2i+1,M}^h \\
\vdots & \vdots & \ddots & \vdots \\
\widehat{x}_{N,0}^h & \widehat{x}_{N,1}^h & \cdots & \widehat{x}_{N,M}^h \ear\in\mathbb{R}^{n_h(N+1) \times (M+1)}.
\eean
This completes the computation of the horizontal states. The same procedure must be applied in the
vertical direction to get the vertical state estimates.
\bibliographystyle{apacite}
\renewcommand\bibliographytypesize{\fontsize{10}{12}\selectfont}
% \bibliography{Bibliography_Ramos}

\begin{thebibliography}{}

\bibitem [\protect \citeauthoryear {%
Ramos%
\ \BBA {} Merc{\`e}re%
}{%
Ramos%
\ \BBA {} Merc{\`e}re%
}{%
{\protect \APACyear {2016a}}%
}]{%
Ramos2016a}
\APACinsertmetastar {%
Ramos2016a}%
\begin{APACrefauthors}%
Ramos, J\BPBI A.%
\BCBT {}\ \BBA {} Merc{\`e}re, G.%
\end{APACrefauthors}%
\unskip\
\newblock
\APACrefYearMonthDay{2016a}{October}{}.
\newblock
{\BBOQ}\APACrefatitle {Image modeling based on a {2-D} stochastic subspace
  system identification algorithm} {Image modeling based on a {2-D} stochastic
  subspace system identification algorithm}.{\BBCQ}
\newblock
\APACjournalVolNumPages{Multidimensional Systems and Signal
  Processing}{28}{}{1133--1165}.
\PrintBackRefs{\CurrentBib}

\bibitem [\protect \citeauthoryear {%
Ramos%
\ \BBA {} Merc\`{e}re%
}{%
Ramos%
\ \BBA {} Merc\`{e}re%
}{%
{\protect \APACyear {2016b}}%
}]{%
Ramos2016b}
\APACinsertmetastar {%
Ramos2016b}%
\begin{APACrefauthors}%
Ramos, J\BPBI A.%
\BCBT {}\ \BBA {} Merc\`{e}re, G.%
\end{APACrefauthors}%
\unskip\
\newblock
\APACrefYearMonthDay{2016b}{}{}.
\newblock
{\BBOQ}\APACrefatitle {Subspace algorithms for identifying
  separable-in-denominator 2D systems with deterministic-stochastic inputs}
  {Subspace algorithms for identifying separable-in-denominator 2d systems with
  deterministic-stochastic inputs}.{\BBCQ}
\newblock
\APACjournalVolNumPages{International Journal of Control}{89}{12}{2584--2610}.
\newblock
\begin{APACrefURL} \url{http://dx.doi.org/10.1080/00207179.2016.1172258}
  \end{APACrefURL}
\PrintBackRefs{\CurrentBib}

\bibitem [\protect \citeauthoryear {%
Ramos%
\ \BBA {} Merc\`{e}re%
}{%
Ramos%
\ \BBA {} Merc\`{e}re%
}{%
{\protect \APACyear {2017a}}%
}]{%
Ramos2017a}
\APACinsertmetastar {%
Ramos2017a}%
\begin{APACrefauthors}%
Ramos, J\BPBI A.%
\BCBT {}\ \BBA {} Merc\`{e}re, G.%
\end{APACrefauthors}%
\unskip\
\newblock
\APACrefYearMonthDay{2017a}{}{}.
\newblock
{\BBOQ}\APACrefatitle {A Stochastic Subspace System Identification Algorithm
  for State Space Systems in the General 2-D Roesser Model Form} {A stochastic
  subspace system identification algorithm for state space systems in the
  general 2-d roesser model form}.{\BBCQ}
\newblock
\APACjournalVolNumPages{International Journal of Control}{}{submitted for
  publication}{}.
\PrintBackRefs{\CurrentBib}

\end{thebibliography}

\end{document}